\documentclass{amsart}
\usepackage[T1]{fontenc}
\usepackage[utf8]{inputenc}
\usepackage{lmodern}
\usepackage[british]{babel}
\usepackage{mathtools,amssymb,amsthm}
\usepackage{array}
\usepackage{needspace}
\usepackage{microtype}
\usepackage{graphicx}
\usepackage{placeins}
\usepackage[hidelinks]{hyperref}
\hypersetup{
  pdftitle={Certified Finite Element Error Bounds for Diffusion Equations with Measure-Valued Principal-Part Delays},
  pdfauthor={Hiroki Ishizaka},
  pdfkeywords={measure-valued principal-part delay, finite element method, a posteriori error estimate, equilibrated flux, Wasserstein distance, delay-law approximation, delay-law certification}
}
\usepackage[nameinlink,capitalize,noabbrev]{cleveref}
\usepackage{etoolbox}

\numberwithin{equation}{section}
\theoremstyle{plain}
\newtheorem{theorem}{Theorem}[section]
\newtheorem{lemma}[theorem]{Lemma}
\newtheorem{proposition}[theorem]{Proposition}
\newtheorem{corollary}[theorem]{Corollary}
\theoremstyle{definition}

\theoremstyle{remark}
\newtheorem{remark}[theorem]{Remark}
\crefname{theorem}{Theorem}{Theorems}
\Crefname{theorem}{Theorem}{Theorems}
\crefname{lemma}{Lemma}{Lemmas}
\Crefname{lemma}{Lemma}{Lemmas}
\crefname{proposition}{Proposition}{Propositions}
\Crefname{proposition}{Proposition}{Propositions}
\crefname{corollary}{Corollary}{Corollaries}
\Crefname{corollary}{Corollary}{Corollaries}
\crefname{definition}{Definition}{Definitions}
\Crefname{definition}{Definition}{Definitions}
\crefname{assumption}{Assumption}{Assumptions}
\Crefname{assumption}{Assumption}{Assumptions}
\crefname{remark}{Remark}{Remarks}
\Crefname{remark}{Remark}{Remarks}
\crefname{example}{Example}{Examples}
\Crefname{example}{Example}{Examples}
\newcommand{\R}{\mathbb{R}}
\newcommand{\M}{\mathcal{M}}
\newcommand{\W}{\mathcal{W}}
\newcommand{\K}{\mathcal{K}}

\newcommand{\dd}{\,\mathrm{d}}
\newcommand{\Tend}{\mathfrak{T}}
\newcommand{\norm}[1]{\left\lVert #1 \right\rVert}
\newcommand{\abs}[1]{\left\lvert #1 \right\rvert}
\newcommand{\pair}[2]{\left\langle #1,#2 \right\rangle}
\DeclareMathOperator{\supp}{supp}

\begin{document}
\raggedbottom
\title[Certified FEM bounds for diffusion with delays]{Certified Finite Element Error Bounds for Diffusion Equations with Measure-Valued Principal-Part Delays}
\author{Hiroki Ishizaka}
\address{Team FEM, Matsuyama, Japan}
\email{h.ishizaka005@gmail.com}
\urladdr{https://teamfem.github.io/hiroki\_ishizaka/}
\begin{abstract}
We derive computable a posteriori finite element error bounds for diffusion equations whose highest-order spatial operator carries a non-negative measure-valued delay. Although the delay law is replaced by finitely many atoms in computations, the bounds measure the error with respect to the solution for the original law. An equilibrated Raviart--Thomas flux controls the spatial residual, and a transport estimate on the computed time reconstruction controls the delay-law residual without additional temporal regularity of the exact solution. A cumulative-mass representation preserves cancellation in this residual and allows a jump between the numerical initial state and the prescribed history. The constants are independent of the mesh aspect ratios and of the number of atoms, and inexact flux equilibration and errors in the elementwise residual means are included. We also give a quantile-based budget for the delay residual and a weighted finite-horizon bound beyond the mass-subcritical regime. Under an elliptic equilibration hypothesis, we prove a two-sided estimate for a combined error norm that includes the time derivative. On right-triangular meshes satisfying the maximum-angle condition, however, the lowest-order Raviart--Thomas estimator can exceed the combined graph error by a factor at least proportional to the aspect ratio. Numerical experiments on anisotropic meshes, including a law beyond the mass-subcritical regime, illustrate the behaviour of the residual indicators and the limitations of efficiency.
\end{abstract}
\maketitle

\section{Introduction}\label{sec:introduction}
Let $\Omega\subset\R^d$, $d\in\{1,2,3\}$, be a bounded Lipschitz domain, and let $V=H_0^1(\Omega)$ and $V^*=H^{-1}(\Omega)$. For a final time $\Tend>0$ and maximal delay $\tau>0$, consider
\begin{align}\label{eq:intro-problem}
  \partial_tu+A_0u+\int_{[0,\tau]}A_1\tilde u(t-s)\,\dd\mu(s)=f,
  \quad 0<t<\Tend,\quad u(0)=u_0.
\end{align}
Here, $A_0$ is coercive and $A_1$ acts through the same energy space. The finite signed Borel measure $\mu$ describes the delay distribution, and $\tilde u$ joins $u$ to a prescribed past history. Delays in the principal part differ from delays in lower-order terms, because high spatial frequencies enter the delayed term; see \cite{YongPan1993,DiBlasioKunischSinestrari1984,KhusainovPokojovyRacke2015}. The a posteriori analysis below concerns non-negative original and atomic delay measures of equal mass, with no atom at the origin. The continuous well-posedness framework also accommodates signed measures.

In practice, $\mu$ is replaced by an atomic law $\mu_N$ before the problem is discretised in space and time. The error that matters for the original problem is $u_\mu-u_h$, not the error relative to the solution for $\mu_N$. Comparing the exact solutions for $\mu$ and $\mu_N$ in a Wasserstein distance may require time regularity of a trajectory that is not known. We instead evaluate the measure perturbation on the numerical time reconstruction, which is available. This gives an a posteriori bound with respect to the original law without additional time regularity of $u_\mu$.

Our first result is a reliable error bound based on splitting the residual into a space--time residual for $\mu_N$ and a delay-law residual. An equilibrated Raviart--Thomas (RT) flux bounds the former with no mesh-shape factor, and a transport estimate bounds the latter in terms of the computed reconstruction; a jump between $u_h(0)$ and the prescribed history is allowed. An exact cumulative-mass identity keeps the cancellations that a transport bound alone would lose. Inexact flux equilibration and errors in the elementwise residual means enter the computable majorant as separate terms. For quantile atoms, the delay-law contribution on a given reconstruction is of order $N^{-1}$; this bounds one term of the majorant and is not a convergence rate for the finite element error.

Our second result concerns efficiency. For grid-aligned delays and scalar constant coefficients, an elliptic equilibration hypothesis yields a two-sided estimate in a norm that controls both the energy error and its time derivative. The hypothesis cannot hold uniformly on anisotropic meshes, even under the maximum-angle condition: on a family of right-triangular meshes, the lowest-order RT estimator exceeds the combined error by a factor at least proportional to the aspect ratio. Reliability on anisotropic meshes therefore does not imply aspect-ratio-uniform efficiency.

A posteriori estimates for fully discrete delay parabolic equations are developed in \cite{WangYiXiao2020,WangYi2022}. In particular, Wang and Yi~\cite{WangYi2022} use delay-dependent elliptic reconstruction to obtain optimal $L^\infty(L^2)$ bounds. Our concern is different: the delay measure used in the computation is itself approximated by atoms, and its residual must be included when estimating the error against the original problem. We evaluate that residual on the computed reconstruction, alongside the space--time residual. A posteriori estimates for the heat equation are developed in \cite{Verfurth2003}, and equilibrated-flux estimates for parabolic problems in \cite{ErnVohralik2010,ErnSmearsVohralik2019}. The difficulties of anisotropic efficiency are discussed in \cite{Kunert2000,Grosman2006,Kopteva2017,Kopteva2020}; in particular, Kunert~\cite{Kunert2000} uses a mesh-alignment matching function, while Kopteva~\cite{Kopteva2020} obtains lower estimates on partially structured anisotropic meshes. Earlier CR--RT error analysis beyond the maximum-angle condition and anisotropic RT interpolation estimates~\cite{IshizakaKobayashiTsuchiya2021,IshizakaRT2022} provide related spatial background. The elliptic CR--RT flux certification and Marini-based estimates in~\cite{IshizakaCRRT2026,IshizakaMariniL2_2026} address different problems; the present bounds also account for the approximation of the delay law.

The continuous well-posedness framework used in this paper was developed in \cite{IshizakaMixed2026}. In particular, we use its well-posedness and energy stability results for equations with measure-valued delays in the principal part. The related preprint \cite{IshizakaWasserstein2026} studies solution-to-solution stability under perturbations of the delay law and small-delay thresholds. Those perturbation estimates are not used here. Instead, we estimate the delay-law residual on the computed reconstruction and apply the energy stability estimate to obtain an a posteriori bound for the original problem.

Sharp spectral sensitivity estimates for measure-valued delays are studied in \cite{Shikhman2026}, while weak-solution theory for nonlinear delayed diffusion is developed in \cite{Tsukamoto2026}. These works address different questions from the present finite element error analysis.

\Cref{sec:framework} fixes the problem, measure notation and transport distance. \Cref{sec:certified-fem} proves reliability and formulates a computable majorant. \Cref{sec:efficiency} establishes the conditional two-sided estimate and the anisotropic obstruction. \Cref{sec:numerics} presents quadrature-based numerical diagnostics and illustrates the limitations of efficiency.

\section{Problem setting and stability margin}\label{sec:framework}
Let $V$ and $H$ be real separable Hilbert spaces with $V\hookrightarrow H$ continuously and densely; identify $H$ with its dual. Suppose that $a_0,a_1:V\times V\to\R$ satisfy
\begin{align}
 |a_i(w,v)|&\le\Lambda_i\|w\|_V\|v\|_V,\quad i\in\{0,1\},\label{eq:form-bounded}\\
 a_0(v,v)&\ge\alpha_0\|v\|_V^2,\label{eq:a0-coercive}
\end{align}
where $\alpha_0,\Lambda_0,\Lambda_1>0$. Let $A_i:V\to V^*$ be the associated operators. We denote by $\M([0,\tau])$ the finite signed Borel measures and we set $\M^0([0,\tau]) := \{\sigma\in\M([0,\tau]):|\sigma|(\{0\})=0\}$ and $\M_+^0 := \M^0\cap\M_+$, where $\M_+$ denotes the non-negative measures. For $\beta\ge0$, we set
\begin{align}
 \|g\|_{L^2_\beta(I;X)}&:=\Bigl(\int_I e^{-2\beta t}\|g(t)\|_X^2\,\dd t\Bigr)^{1/2},\label{eq:weighted-norm}\\
 \kappa_\sigma(\beta)&:=\int_{[0,\tau]}e^{-\beta s}\,\dd|\sigma|(s),\quad
 d_{\mu,\beta}:=\alpha_0-\Lambda_1\kappa_\mu(\beta).\label{eq:kappa}
\end{align}
For $u\in L^2(0,\Tend;V)$ and a prescribed $\psi\in L^2(-\tau,0;V)$, let
\begin{align}\label{eq:joined-trajectory}
 \tilde u(t):=\begin{cases}\psi(t),&-\tau<t<0,\\u(t),&0<t<\Tend.\end{cases}
\end{align}
We write $\W(0,\Tend):=\{u\in L^2(0,\Tend;V):\partial_tu\in L^2(0,\Tend;V^*)\}$ with its natural graph norm. It is continuously embedded in $\mathcal C([0,\Tend];H)$. For $\sigma\in\M([0,\tau])$, the history operator
\begin{align}\label{eq:history-operator}
 (\K_\sigma g)(t):=\int_{[0,\tau]} A_1g(t-s)\,\dd\sigma(s)
\end{align}
is well defined in $L^2(0,\Tend;V^*)$ for $g\in L^2(-\tau,\Tend;V)$, independently of the strongly measurable representative, and satisfies
\begin{align}\label{eq:measure-convolution-bound}
 \|\K_\sigma g\|_{L^2(0,\Tend;V^*)}\le\Lambda_1\|\sigma\|_{\M([0,\tau])}\|g\|_{L^2(-\tau,\Tend;V)};
\end{align}
see \cite[Lemma~2.3]{IshizakaMixed2026}. We study
\begin{align}\label{eq:problem}
 \partial_tu+A_0u+\K_\mu\tilde u=f\quad\text{in }V^*,\quad u(0)=u_0\ \text{in }H,
\end{align}
where $\mu\in\M^0([0,\tau])$, $u_0\in H$, $f\in L^2(0,\Tend;V^*)$ and $\psi\in L^2(-\tau,0;V)$. For each such fixed $\mu$, the problem has a unique solution in $\W(0,\Tend)$ by \cite[Theorem~3.6]{IshizakaMixed2026}. If a measure has an atom $m_0\delta_0$ at the origin, this term is absorbed into $A_0+m_0A_1$, provided coercivity is retained. In the non-negative case, $\kappa_\mu(0)=M:=\mu([0,\tau])$ and $d_{\mu,0}=\alpha_0-\Lambda_1M$. We call such a law \emph{mass-subcritical} if $\Lambda_1M<\alpha_0$; then $\beta=0$ is admissible and $d_{\mu,0}$ does not depend on the delay locations. If $\Lambda_1M\ge\alpha_0$, a positive margin requires $\beta>0$. More generally, $|\mu|(\{0\})=0$ implies $\kappa_\mu(\beta)\to0$ as $\beta\to\infty$ by dominated convergence; thus $d_{\mu,\beta}>0$ for some $\beta\ge0$. The energy argument in \Cref{thm:fem-certificate} needs only such a positive margin, not a delay-uniform stability estimate.

\subsection{Diffusion realisation}\label{sec:pde-realisation}
For the finite element analysis, let $\Omega\subset\R^d$, $d\in\{1,2,3\}$, be bounded and Lipschitz, and we set $H := L^2(\Omega)$, $V := H_0^1(\Omega)$, $V^* := H^{-1}(\Omega)$ with $\|v\|_V=\|\nabla v\|_{L^2(\Omega;\R^d)}$. Suppose that $B_0,B_1\in L^\infty(\Omega;\R^{d\times d})$ and that $\xi^\top B_0(x)\xi\ge\lambda_0|\xi|^2$ for almost every $x$ and any $\xi\in\R^d$. We define
\begin{align}\label{eq:diffusion-forms}
 a_i(w,v):=\int_\Omega B_i(x)\nabla w\cdot\nabla v\,\dd x,\quad i\in\{0,1\}.
\end{align}
Then, \eqref{eq:form-bounded}--\eqref{eq:a0-coercive} hold with $\alpha_0=\lambda_0$ and $\Lambda_i=\|B_i\|_{L^\infty(\Omega;\mathcal L(\R^d))}$; see \cite[Proposition~4.1]{IshizakaMixed2026}.

\subsection{Transport notation}
For $\mu,\nu\in\M_+([0,\tau])$ with the same mass $M>0$, let $\Pi(\mu,\nu)$ be the non-negative Borel measures on $[0,\tau]^2$ whose marginals are $\mu$ and $\nu$. We define the equal-mass transport cost as
\begin{align}\label{eq:W1}
 W_1(\mu,\nu):=\inf_{\gamma\in\Pi(\mu,\nu)}\int_{[0,\tau]^2}|s-r|\,\dd\gamma(s,r).
\end{align}
In one dimension, if $F_\mu(t)=\mu([0,t])$ and $F_\nu(t)=\nu([0,t])$, then
\begin{align*}
W_1(\mu,\nu)
=\int_0^\tau |F_\mu(t)-F_\nu(t)|\,\dd t;
\end{align*}
see \cite{Santambrogio2015}. We write
\begin{align}\label{eq:cdf-discrepancy}
 \Delta_{\nu,\mu}(t)
 :=(\nu-\mu)([0,\min\{t,\tau\}]),
 \quad t\ge0,
\end{align}
for the cumulative discrepancy. It vanishes for $t\ge\tau$, and the identity above reads $W_1(\mu,\nu)=\int_0^\tau|\Delta_{\nu,\mu}(t)|\,\dd t$.

\section{Finite element error bounds}\label{sec:certified-fem}
Throughout, the error is measured against the solution for the \emph{original} delay law. The residual splits into the residual $R_N$ for the atomic law and a delay-law term that involves only the computed reconstruction:
\begin{align*}
 R=R_N+\K_{\mu_N-\mu}\tilde u_h.
\end{align*}
The equilibrated flux controls $R_N$; the second term is bounded by transport plans and, more sharply, by a cumulative-mass identity. Only the stability margin $d_{\mu,\beta}$ enters, and no regularity estimate for the exact solution is used.

\subsection{Reliability with equilibrated fluxes}
\label{subsec:fem-reliability}
For the finite element analysis, $\Omega\subset\R^d$, $d\in\{1,2,3\}$, is a bounded polyhedral Lipschitz domain, not necessarily convex. Let $\mathbb{T}_h$ be a simplicial mesh of $\overline{\Omega}$ composed of closed $d$-simplices. Let $N_e := \#\mathbb{T}_h$ denote the number of elements of the mesh and write
\begin{align*}
\displaystyle
\mathbb{T}_h = \{ T_j \}_{j=1}^{N_e}, \quad \overline{\Omega} = \bigcup_{j=1}^{N_e} T_j.
\end{align*}
For each $j \in \{1,\ldots,N_e\}$, we set
\begin{align*}
\displaystyle
h_{T_j} := \operatorname{diam}(T_j), \quad h := \max_{1 \le j \le N_e} h_{T_j}.
\end{align*}
For simplicity, we assume that $\mathbb{T}_h$ is conforming; that is, $\mathbb{T}_h$ is a simplicial mesh of $\overline{\Omega}$ with no hanging nodes. Throughout this paper, we assume $0 < h \leq 1$. For any $E \subset \mathbb R^d$, we write $|E|_s :=\mathcal H^s(E)$, where $\mathcal H^s$ is the $s$-dimensional Hausdorff measure. Let $\mathbb P^k(T)$ be the polynomials of total degree at most $k$. We write $\mathbb P^k(\mathbb T_h)$ for the corresponding discontinuous piecewise polynomials and $\Pi_h^0$ for the $L^2$-orthogonal projection onto $\mathbb P^0(\mathbb T_h)$, and let $V_h\subset H_0^1(\Omega)$ be a conforming finite element space, for example $\mathbb P^1(\mathbb T_h)\cap H_0^1(\Omega)$. With the centroid $x_T$ of $T$, the local lowest-order Raviart--Thomas space $\mathbb{RT}^{0}(T)$ consists of the fields $\boldsymbol a_T+c_T(x-x_T)$, $\boldsymbol a_T\in\R^d$, $c_T\in\R$, with divergence $dc_T$. The global space $\mathbb{RT}^{0}(\mathbb T_h)$ consists of the fields in $H(\operatorname{div};\Omega)$ whose restriction to each $T$ lies in $\mathbb{RT}^{0}(T)$, and no normal flux is prescribed on $\partial\Omega$. The reliability estimates assume neither shape regularity, $h_T^d/|T|_d\le\gamma_0$, nor, for $d\in\{2,3\}$, the semi-regular condition $H_T/h_T\le\gamma_1$ of \cite{IshizakaKobayashiTsuchiya2023}, where $H_T:=h_T\prod_{i=1}^dh_{T,i}/|T|_d$ with the characteristic edge lengths $h_{T,i}$ of that reference, nor a maximum-angle or bounded-aspect-ratio condition; increasingly anisotropic elements are allowed.

Suppose that a numerical method has produced a reconstruction
\begin{align}\label{eq:fem-reconstruction}
  u_h\in H^1(0,\Tend;V),\quad u_h(t)\in V_h\quad\text{for a.e. }t.
\end{align}
A continuous, piecewise affine time interpolant of conforming finite element solutions, such as those obtained by the backward Euler method, is one example. On $(-\tau,0)$, we set $\tilde u_h=\psi$, where the prescribed history satisfies $\psi\in H^1(-\tau,0;V)$.  Thus, $\tilde u_h$ is piecewise $H^1$ in time, but it need not be continuous at $t=0$.  We set
\begin{align}\label{eq:join-jump}
  J_h:=u_h(0)-\psi(0)\in V,\quad
  D_h:=\left(\norm{\partial_t\psi}_{L^2(-\tau,0;V)}^2
       +\norm{\partial_t u_h}_{L^2(0,\Tend;V)}^2\right)^{1/2}.
\end{align}
The traces in \eqref{eq:join-jump} are well defined in $V$.  In particular, a non-zero $J_h$ is allowed even if the exact initial datum is compatible with its history.

Let $\mu_N\in\M_+^0([0,\tau])$ be an atomic measure of the same mass $M>0$ as $\mu$. We define the numerical diffusive flux and the volumetric residual datum as
\begin{align}\label{eq:numerical-flux}
  F_N(t)&:=B_0\nabla u_h(t)
  +\int_{[0,\tau]}B_1\nabla\tilde u_h(t-s)\,\dd\mu_N(s),\\
  g(t)&:=\partial_t u_h(t)-f(t),\quad g_h(t):=\Pi_h^0g(t).
\end{align}
When $f\in L^2(0,\Tend;L^2(\Omega))$, these expressions belong to the required $L^2$ spaces. We choose an equilibrated flux
\begin{align}\label{eq:equilibrated-flux}
  \sigma_h\in L^2(0,\Tend;\mathbb{RT}^{0}({\mathbb T}_h)),\quad
  \nabla\cdot\sigma_h(t)=g_h(t)\quad\text{for a.e. }t.
\end{align}
Such a flux exists because the discrete divergence maps the unrestricted global $\mathbb{RT}^{0}$ space onto the elementwise constants.  It may be computed by a mixed or constrained least-squares solve; no particular construction is needed for reliability.

For $\beta\ge0$, we introduce
\begin{align}\label{eq:fem-estimators}
  \eta_{\mathrm{flux}}&:=
  \norm{F_N-\sigma_h}_{L^2_\beta(0,\Tend;L^2(\Omega;\R^d))},\\
  \eta_{\mathrm{osc}}&:=\left(
    \int_0^{\Tend}e^{-2\beta t}
      \sum_{T\in{\mathbb T}_h}\frac{h_T^2}{\pi^2}
      \norm{g(t)-g_h(t)}_{L^2(T)}^2\,\dd t
  \right)^{1/2},\label{eq:eta-osc}\\
  \eta_{\mathrm{ker}}&:=\Lambda_1\left[
    W_1(\mu,\mu_N)D_h
    +\sqrt{M W_1(\mu,\mu_N)}\,\norm{J_h}_V\right],\label{eq:eta-ker}\\
  \eta_{\mathrm{init}}&:=\norm{u_0-u_h(0)}_H.\label{eq:eta-init}
\end{align}
The constant $h_T/\pi$ follows from the convex-domain Poincar\'e inequality~\cite{PayneWeinberger1960,Bebendorf2003} and does not depend on the aspect ratio.

\begin{theorem}[Reliability with a measure-valued principal-part delay]
\label[theorem]{thm:fem-certificate}
Assume the diffusion setting of \Cref{sec:pde-realisation}, with $f\in L^2(0,\Tend;L^2(\Omega))$, $u_0\in H$, and $\psi\in H^1(-\tau,0;V)$.  Let $\mu,\mu_N\in\M_+^0([0,\tau])$ have the same mass $M>0$.  Suppose that \eqref{eq:fem-reconstruction}--\eqref{eq:equilibrated-flux} hold.  Choose $\beta\ge0$ such that
\begin{align}\label{eq:fem-stability-margin}
  d_{\mu,\beta}=\alpha_0-\Lambda_1\kappa_\mu(\beta)>0.
\end{align}
Then, the exact solution $u_\mu$ of \eqref{eq:problem} and the reconstruction $u_h$ satisfy
\begin{align}\label{eq:fem-reliability}
  \norm{u_\mu-u_h}_{L^2_\beta(0,\Tend;V)}
  \le \frac{\eta_{\mathrm{init}}}{\sqrt{2d_{\mu,\beta}}}
   +\frac{\eta_{\mathrm{flux}}+\eta_{\mathrm{osc}}+\eta_{\mathrm{ker}}}
          {d_{\mu,\beta}}.
\end{align}
All constants are explicit, and the reliability constant is independent of the spatial mesh and of the number and locations of the delay atoms. The right-hand side involves only the data, the reconstruction, the flux and an upper bound for the transport cost.
\end{theorem}

\begin{proof}
We set $e:=u_\mu-u_h$ for $t>0$ and extend $e$ by zero to $(-\tau,0)$, because $u_\mu$ and $u_h$ use the same prescribed history there.  In $L^2(0,\Tend;V^*)$, the error satisfies
\begin{align}\label{eq:fem-error-equation}
  \partial_te+A_0e+\K_\mu\tilde e&=R,\\
  R&:=f-\partial_t u_h-A_0u_h-\K_\mu\tilde u_h.
\end{align}
The initial error is $e(0)=u_0-u_h(0)$ in $H$. We write
\begin{align*}
  R=\underbrace{f-\partial_t u_h-A_0u_h-\K_{\mu_N}\tilde u_h}_{R_N}
      +\K_{\mu_N-\mu}\tilde u_h.
\end{align*}
By the diffusion-form identity and \eqref{eq:numerical-flux}, for $v\in V$ and almost every $t$,
\begin{align*}
  \pair{R_N(t)}{v}
  &=-(g(t),v)_{L^2(\Omega)}-(F_N(t),\nabla v)_{L^2(\Omega;\R^d)}.
\end{align*}
Integration by parts in \eqref{eq:equilibrated-flux} gives
\begin{align*}
  \pair{R_N(t)}{v}
  &=-(g(t)-g_h(t),v)_{L^2(\Omega)}
    -(F_N(t)-\sigma_h(t),\nabla v)_{L^2(\Omega;\R^d)}.
\end{align*}
Because $g-g_h$ has zero mean on any simplex $T$, the convex-element Poincar\'e inequality yields
\begin{align*}
  \abs{(g-g_h,v)_{L^2(\Omega)}}
  &\le \sum_{T\in{\mathbb T}_h}
    \norm{g-g_h}_{L^2(T)}\norm{v-v_T}_{L^2(T)}\\
  &\le\left(\sum_{T\in{\mathbb T}_h}
    \frac{h_T^2}{\pi^2}\norm{g-g_h}_{L^2(T)}^2\right)^{1/2}\norm{v}_V,
\end{align*}
where $v_T$ is the mean of $v$ on $T$.  Consequently,
\begin{align}\label{eq:fem-residual-flux}
  \norm{R_N}_{L^2_\beta(0,\Tend;V^*)}
  \le\eta_{\mathrm{flux}}+\eta_{\mathrm{osc}}.
\end{align}

For the kernel term, we fix a transport plan $\gamma\in\Pi(\mu_N,\mu)$.  The equality of the marginals implies
\begin{align*}
  \K_{\mu_N-\mu}\tilde u_h(t)
  =\int_{[0,\tau]^2}A_1
       \bigl(\tilde u_h(t-s)-\tilde u_h(t-r)\bigr)\,\dd\gamma(s,r).
\end{align*}
With $J_h=u_h(0)-\psi(0)$ from \eqref{eq:join-jump}, we define $w_h(t):=\tilde u_h(t)-J_h\mathbf{1}_{[0,\infty)}(t)$. Then, $w_h\in H^1(-\tau,\Tend;V)$ and
\begin{align*}
\norm{\partial_t w_h}_{L^2(-\tau,\Tend;V)}=D_h.
\end{align*}
For $0\le s\le r\le\tau$, the fundamental theorem of calculus and Minkowski's inequality give
\begin{align*}
\norm{w_h(\cdot-s)-w_h(\cdot-r)}
       _{L^2(0,\Tend;V)}
\le
\int_s^r \norm{\partial_t w_h(\cdot-\xi)}
       _{L^2(0,\Tend;V)}\,\dd\xi
\le (r-s)D_h.
\end{align*}
The jump contribution is supported on $[s,r)\cap(0,\Tend)$ and therefore satisfies
\begin{align*}
\norm{J_h\left(
\mathbf{1}_{[s,\infty)}
-\mathbf{1}_{[r,\infty)}
\right)}_{L^2(0,\Tend;V)}
\le (r-s)^{1/2}\norm{J_h}_V.
\end{align*}
Interchanging $s$ and $r$ if necessary, we obtain
\begin{align}\label{eq:jump-transport}
\norm{\tilde u_h(\cdot-s)-\tilde u_h(\cdot-r)}
     _{L^2(0,\Tend;V)}
\le
\abs{s-r}D_h
+\abs{s-r}^{1/2}\norm{J_h}_V.
\end{align}
Thus, the exponent $1/2$ arises from the $L^2$ translation estimate for the jump. We write
\begin{align*}
Q_\gamma:=
\int_{[0,\tau]^2}\abs{s-r}\,\dd\gamma(s,r).
\end{align*}
By the boundedness of $A_1$, the fact that the weighted norm does not exceed the unweighted norm, and Minkowski's inequality, we have
\begin{align*}
\norm{\K_{\mu_N-\mu}\tilde u_h}
       _{L^2_\beta(0,\Tend;V^*)}
\le
\Lambda_1\int_{[0,\tau]^2}
\left(
\abs{s-r}D_h
+\abs{s-r}^{1/2}\norm{J_h}_V
\right)\dd\gamma(s,r).
\end{align*}
Because $\gamma([0,\tau]^2)=M$, the Cauchy--Schwarz inequality yields
\begin{align*}
\int_{[0,\tau]^2}\abs{s-r}^{1/2}\,\dd\gamma(s,r)
\le M^{1/2}Q_\gamma^{1/2}.
\end{align*}
Consequently,
\begin{align*}
\norm{\K_{\mu_N-\mu}\tilde u_h}
     _{L^2_\beta(0,\Tend;V^*)}
\le
\Lambda_1\left(
D_hQ_\gamma+
M^{1/2}\norm{J_h}_VQ_\gamma^{1/2}
\right).
\end{align*}
Taking an optimal transport plan on $[0,\tau]$ and using $Q_\gamma=W_1(\mu_N,\mu)$, we conclude that
\begin{align}\label{eq:fem-residual-kernel}
\norm{\K_{\mu_N-\mu}\tilde u_h}
     _{L^2_\beta(0,\Tend;V^*)}
\le\eta_{\mathrm{ker}}.
\end{align}

For the zero-history extension of $e$, a change of variables gives, for any $s\in[0,\tau]$,
\begin{align*}
\norm{\tilde e(\cdot-s)}_{L^2_\beta(0,\Tend;V)}
\le e^{-\beta s}
\norm{e}_{L^2_\beta(0,\Tend;V)}.
\end{align*}
Therefore, Minkowski's integral inequality and the boundedness of $A_1$ yield
\begin{align*}
\norm{\K_\mu\tilde e}_{L^2_\beta(0,\Tend;V^*)}
\le \Lambda_1
\int_{[0,\tau]}e^{-\beta s}\,\dd|\mu|(s)
\norm{e}_{L^2_\beta(0,\Tend;V)}
=\Lambda_1\kappa_\mu(\beta)
\norm{e}_{L^2_\beta(0,\Tend;V)}.
\end{align*}
Testing \eqref{eq:fem-error-equation} with $e^{-2\beta t}e(t)$ and integrating in time,
the Gelfand-triple identity, coercivity of $a_0$, and the preceding bound for the memory term give
\begin{align*}
  d_{\mu,\beta}
  \norm{e}_{L^2_\beta(0,\Tend;V)}^2
  \le \tfrac12\norm{e(0)}_H^2
  +\norm{R}_{L^2_\beta(0,\Tend;V^*)}
   \norm{e}_{L^2_\beta(0,\Tend;V)}.
\end{align*}
Solving this quadratic inequality yields
\begin{align*}
  \norm{e}_{L^2_\beta(0,\Tend;V)}
  \le
  \frac{
    \norm{R}_{L^2_\beta(0,\Tend;V^*)}
    +
    \sqrt{
      \norm{R}_{L^2_\beta(0,\Tend;V^*)}^2
      +2d_{\mu,\beta}\norm{e(0)}_H^2
    }
  }{2d_{\mu,\beta}}.
\end{align*}
Using $\sqrt{a+b}\le\sqrt a+\sqrt b$ for $a,b\ge0$, we obtain
\begin{align*}
  \norm{e}_{L^2_\beta(0,\Tend;V)}
  \le
  \frac{\norm{e(0)}_H}{\sqrt{2d_{\mu,\beta}}}
  +
  \frac{\norm{R}_{L^2_\beta(0,\Tend;V^*)}}
       {d_{\mu,\beta}}.
\end{align*}
Finally, \eqref{eq:fem-residual-flux} and \eqref{eq:fem-residual-kernel}, together with the residual decomposition, complete the proof.
\end{proof}

\subsection{Sharper bounds for the delay-law residual}
\label{subsec:delay-residual}
\begin{proposition}[A cumulative-mass bound for the history jump]
\label[proposition]{prop:cdf-jump}
Under the assumptions of \Cref{thm:fem-certificate}, the estimate \eqref{eq:fem-reliability} holds with $\eta_{\mathrm{ker}}$ replaced by
\begin{align}\label{eq:cdf-kernel-estimator}
 \eta_{\mathrm{ker}}^{\mathrm{cdf}}
 :=\Lambda_1\bigl[W_1(\mu,\mu_N)D_h
       +\|\Delta_{\mu_N,\mu}\|_{L^2_\beta(0,\Tend)}\|J_h\|_V\bigr],
\end{align}
where $\Delta_{\mu_N,\mu}$ is the cumulative discrepancy \eqref{eq:cdf-discrepancy}. Furthermore,
\begin{align}\label{eq:cdf-W1-control}
 \|\Delta_{\mu_N,\mu}\|_{L^2_\beta(0,\Tend)}\le
 \sqrt{M W_1(\mu,\mu_N)},
\end{align}
so that this estimator never exceeds $\eta_{\mathrm{ker}}$.
\end{proposition}

\begin{proof}
We write $\tilde u_h=A+J_h\mathbf1_{[0,\infty)}$ on $(-\tau,\Tend)$, with $A\in H^1(-\tau,\Tend;V)$ and $\|A'\|_{L^2(-\tau,\Tend;V)}=D_h$. For a coupling $\gamma\in\Pi(\mu_N,\mu)$, the integral of the continuous part $A(t-s)-A(t-r)$ is bounded in $L^2_\beta(0,\Tend;V)$ by $D_h\int|s-r|\,\dd\gamma$. The step part is exactly $J_h\Delta_{\mu_N,\mu}(t)$ for almost every $t>0$. Taking an optimal coupling gives \eqref{eq:cdf-kernel-estimator}. On the line, the transport identity for equal-mass measures gives $W_1(\mu,\mu_N)=\int_0^\tau|\Delta_{\mu_N,\mu}(t)|\,\dd t$.  Because the measures are non-negative and have the same mass $M$, we have $|\Delta_{\mu_N,\mu}(t)|\le M$. Furthermore, $\Delta_{\mu_N,\mu}(t)=0$ for $t\ge\tau$. Therefore, $\beta\ge0$ and the one-dimensional transport identity give
\begin{align*}
 \|\Delta_{\mu_N,\mu}\|_{L^2_\beta(0,\Tend)}^2
 \le M\int_0^\tau
       |\Delta_{\mu_N,\mu}(t)|\,\dd t
 =M W_1(\mu,\mu_N),
\end{align*}
which proves \eqref{eq:cdf-W1-control}.
\end{proof}

\begin{proposition}[Cancellation-preserving delay-law residual]
\label[proposition]{prop:cancellation-residual}
Under the assumptions of \Cref{thm:fem-certificate}, we set $\delta:=\mu_N-\mu$ and $\Delta:=\Delta_{\mu_N,\mu}$. Let $\dot{\tilde u}_h$ denote the piecewise time derivative, equal to $\partial_t\psi$ on $(-\tau,0)$ and $\partial_tu_h$ on $(0,\Tend)$. We define the known $V$-valued function
\begin{align}\label{eq:cancellation-reconstruction}
 \mathcal C_h(t)
 :=\int_0^\tau
   \Delta(s)\dot{\tilde u}_h(t-s)\,\dd s
   +J_h\Delta(t),
 \quad 0<t<\Tend.
\end{align}
Then, for almost every $t\in(0,\Tend)$,
\begin{align}\label{eq:cancellation-identity}
 \K_{\mu_N-\mu}\tilde u_h(t)
 =A_1\mathcal C_h(t)
 \quad\text{in }V^*.
\end{align}
Furthermore, the estimator
\begin{align}\label{eq:eta-kernel-rec}
 \eta_{\mathrm{ker}}^{\mathrm{rec}}(\beta)
 :=\Lambda_1
   \|\mathcal C_h\|_{L^2_\beta(0,\Tend;V)}
 \le\eta_{\mathrm{ker}}^{\mathrm{cdf}}
\end{align}
may replace $\eta_{\mathrm{ker}}^{\mathrm{cdf}}$ in the preceding reliability estimates and in \eqref{eq:computable-majorant}. It depends only on the cumulative discrepancy and the time reconstruction. If $A_1=\Lambda_1 R_V$, where $R_V:V\to V^*$ is the Riesz map, then
\begin{align}\label{eq:eta-rec-exact}
 \eta_{\mathrm{ker}}^{\mathrm{rec}}(\beta)
 =
 \|\K_{\mu_N-\mu}\tilde u_h\|
 _{L^2_\beta(0,\Tend;V^*)}.
\end{align}
\end{proposition}

\begin{proof}
Because $\mu_N$ and $\mu$ have the same mass, $\delta([0,\tau])=0$.

We first separate the continuous and jump parts of the time reconstruction. We define
\begin{align*}
 A(t):=\tilde u_h(t)
       -J_h\mathbf1_{[0,\infty)}(t),
 \quad -\tau<t<\Tend.
\end{align*}
The definition of $J_h$ gives
\begin{align*}
 A(0-)=\psi(0)=u_h(0)-J_h=A(0+).
\end{align*}
Consequently,
\begin{align*}
 A\in H^1(-\tau,\Tend;V),\quad
 A'=\dot{\tilde u}_h
 \quad\text{a.e.},
\end{align*}
and
\begin{align*}
 \|A'\|_{L^2(-\tau,\Tend;V)}=D_h.
\end{align*}
We fix $t\in(0,\Tend)$ for which the following integrals are defined. Because $\delta$ has total mass zero,
\begin{align*}
 \int_{[0,\tau]}A(t-s)\,\dd\delta(s)
 =
 \int_{[0,\tau]}
 \left(A(t-s)-A(t-\tau)\right)
 \,\dd\delta(s).
\end{align*}
The absolute continuity of $A$ gives
\begin{align*}
 A(t-s)-A(t-\tau)
 =\int_s^\tau A'(t-r)\,\dd r.
\end{align*}
Thus, Fubini's theorem yields
\begin{align*}
\int_{[0,\tau]}A(t-s)\,\dd\delta(s)
&=
 \int_{[0,\tau]}
 \int_s^\tau A'(t-r)\,\dd r\,\dd\delta(s)
=
 \int_0^\tau
 \delta([0,r])A'(t-r)\,\dd r\\
&=
 \int_0^\tau
 \Delta(r)\dot{\tilde u}_h(t-r)\,\dd r.
\end{align*}
For the jump part, we have
\begin{align*}
 \int_{[0,\tau]}
 \mathbf1_{[0,\infty)}(t-s)\,\dd\delta(s) =
 \delta([0,\min\{t,\tau\}])
 =\Delta(t).
\end{align*}
Here, $\Delta(t)=0$ for $t\ge\tau$ because $\delta([0,\tau])=0$. Combining the two parts, we obtain
\begin{align*}
 \int_{[0,\tau]}
 \tilde u_h(t-s)\,\dd\delta(s) =
 \int_0^\tau
 \Delta(s)\dot{\tilde u}_h(t-s)\,\dd s
 +J_h\Delta(t)
 =\mathcal C_h(t).
\end{align*}
Applying the bounded linear operator $A_1:V\to V^*$ proves \eqref{eq:cancellation-identity}.

It remains to estimate $\mathcal C_h$. For any $s\in[0,\tau]$, a change of variables and $\beta\ge0$ give
\begin{align*}
 \|\dot{\tilde u}_h(\cdot-s)\|
       _{L^2_\beta(0,\Tend;V)} \le
 \|\dot{\tilde u}_h\|
       _{L^2(-\tau,\Tend;V)}
 =D_h.
\end{align*}
Minkowski's inequality therefore gives
\begin{align*}
 \|\mathcal C_h\|_{L^2_\beta(0,\Tend;V)}
 &\le
 \int_0^\tau|\Delta(s)|
 \|\dot{\tilde u}_h(\cdot-s)\|
       _{L^2_\beta(0,\Tend;V)}
 \,\dd s +
 \|J_h\|_V
 \|\Delta\|_{L^2_\beta(0,\Tend)}\\
 &\le
 D_h\int_0^\tau|\Delta(s)|\,\dd s
 +\|J_h\|_V
 \|\Delta\|_{L^2_\beta(0,\Tend)}.
\end{align*}
By the one-dimensional transport identity,
\begin{align*}
 \int_0^\tau|\Delta(s)|\,\dd s
 =W_1(\mu,\mu_N).
\end{align*}
Therefore,
\begin{align*}
 \Lambda_1
 \|\mathcal C_h\|_{L^2_\beta(0,\Tend;V)}
 \le\eta_{\mathrm{ker}}^{\mathrm{cdf}}.
\end{align*}
Together with \eqref{eq:cancellation-identity} and the boundedness of $A_1$, this proves \eqref{eq:eta-kernel-rec} and the asserted reliability bounds.

Finally, if $A_1=\Lambda_1R_V$, the isometry property of the Riesz map gives
\begin{align*}
 \|A_1v\|_{V^*}
 =\Lambda_1\|v\|_V
 \quad \forall v\in V.
\end{align*}
Applying this identity to $v=\mathcal C_h(t)$ and integrating in time proves \eqref{eq:eta-rec-exact}.
\end{proof}

\begin{remark}[An obstruction to efficiency of transport bounds]
\label[remark]{rem:transport-not-efficient}
The transport-based estimator may be strictly positive even when the delay-law residual vanishes. Let $0<a<b<\tau$ and set
\begin{align*}
 \mu=\frac M2(\delta_a+\delta_b),
 \quad
 \mu_N=M\delta_{(a+b)/2}.
\end{align*}
For the compatible reconstruction $\tilde u_h(t)=tv$, with $0\ne v\in V_h$, the two measures have the same mass and first moment. Consequently,
\begin{align*}
 \K_{\mu_N-\mu}\tilde u_h=0,
 \quad
 \eta_{\mathrm{ker}}^{\mathrm{rec}}=0.
\end{align*}
However,
\begin{align*}
 W_1(\mu,\mu_N)=\frac{M(b-a)}2>0,
\end{align*}
and $D_h>0$. Thus, $\eta_{\mathrm{ker}}^{\mathrm{cdf}}>0$. Therefore, the transport-based estimator admits no unconditional efficiency bound in terms of the delay-law residual; reliability is not affected.
\end{remark}

\subsection{Inexact equilibration and quadrature-based means}
\label{subsec:inexact-flux}
In practice, the cell means $g_h=\Pi_h^0g$ are assembled by quadrature, and the equilibrium condition in \eqref{eq:equilibrated-flux} may fail when the algebraic system is solved inexactly. Both defects have to enter the error bound. Let $C_{P,\Omega}>0$ be any known upper bound in the Dirichlet Poincar\'e inequality
\begin{align}\label{eq:dirichlet-cp}
 \|v\|_{L^2(\Omega)}\le C_{P,\Omega}\|v\|_V,
 \quad v\in V.
\end{align}
For an arbitrary flux $\sigma_h\in L^2(0,\Tend;\mathbb{RT}^{0}({\mathbb T}_h))$ and any piecewise constant approximation $\hat g_h\in L^2(0,\Tend;\mathbb{P}^{0}({\mathbb T}_h))$ to $g_h$, we define
\begin{align}
 \eta_{\mathrm{flux}}(\sigma_h)&:=
 \|F_N-\sigma_h\|_{L^2_\beta(0,\Tend;L^2(\Omega;\R^d))},\label{eq:eta-flux-sigma}\\
 \eta_{\mathrm{mean}}&:=C_{P,\Omega}
  \|\Pi_h^0g-\hat g_h\|_{L^2_\beta(0,\Tend;L^2(\Omega))},\label{eq:mean-quadrature-defects}\\
 \hat\eta_{\mathrm{div}}&:=C_{P,\Omega}
  \|\hat g_h-\operatorname{div}\sigma_h\|_{L^2_\beta(0,\Tend;L^2(\Omega))}.\label{eq:hat-eta-div}
\end{align}

\begin{proposition}[A computable majorant allowing inexact equilibration]
\label[proposition]{prop:computable-majorant}
We retain the assumptions of \Cref{thm:fem-certificate}, but replace \eqref{eq:equilibrated-flux} by the requirement
\begin{align*}
 \sigma_h \in L^2(0,\Tend;\mathbb{RT}^{0}(\mathbb T_h)),
\end{align*}
without imposing $\operatorname{div}\sigma_h=g_h$. Let $\hat g_h \in L^2(0,\Tend;\mathbb P^0(\mathbb T_h))$ be arbitrary. Then,
\begin{align}\label{eq:computable-majorant}
 \|u_\mu-u_h\|_{L^2_\beta(0,\Tend;V)}
 &\le \frac{\eta_{\mathrm{init}}}{\sqrt{2d_{\mu,\beta}}}
 +\frac{\eta_{\mathrm{flux}}(\sigma_h)+\eta_{\mathrm{osc}}
 +\eta_{\mathrm{mean}}+\hat\eta_{\mathrm{div}}
 +\eta_{\mathrm{ker}}^{\mathrm{cdf}}}{d_{\mu,\beta}}.
\end{align}
For $\hat g_h=g_h$, one has $\eta_{\mathrm{mean}}=0$ and $\hat\eta_{\mathrm{div}}$ is the algebraic equilibration defect; if, in addition, \eqref{eq:equilibrated-flux} holds, then \eqref{eq:computable-majorant} is \eqref{eq:fem-reliability} with $\eta_{\mathrm{ker}}$ replaced by $\eta_{\mathrm{ker}}^{\mathrm{cdf}}\le\eta_{\mathrm{ker}}$. The estimate requires neither a uniform mesh-shape condition nor an exact solution of the algebraic saddle-point system.
\end{proposition}

\begin{proof}
For $v\in V$ and almost every $t\in(0,\Tend)$, integration by parts gives
\begin{align*}
 \langle R_N(t),v\rangle
 &= -(g(t)-g_h(t),v)_{L^2(\Omega)} -(F_N(t)-\sigma_h(t),\nabla v)
    _{L^2(\Omega;\R^d)}\\
&\quad -(g_h(t)-\hat g_h(t),v)_{L^2(\Omega)} -(\hat g_h(t)
       -\operatorname{div}\sigma_h(t),v)
    _{L^2(\Omega)}.
\end{align*}
Because $g_h=\Pi_h^0g$, the function $g-g_h$ has zero mean on each element. The elementwise Poincar\'e inequality therefore gives
\begin{align*}
 |(g(t)-g_h(t),v)_{L^2(\Omega)}|
 \le \omega_{\mathrm{osc}}(t)\|v\|_V,
\end{align*}
where
\begin{align*}
 \omega_{\mathrm{osc}}(t):=
 \left(
 \sum_{K\in\mathbb T_h}
 \frac{h_K^2}{\pi^2}
 \|g(t)-g_h(t)\|_{L^2(K)}^2
 \right)^{1/2}.
\end{align*}
The remaining terms are bounded by Cauchy--Schwarz and \eqref{eq:dirichlet-cp}.
Consequently, 
\begin{align*}
 \|R_N(t)\|_{V^*}
 &\le
 \|F_N(t)-\sigma_h(t)\|_{L^2(\Omega;\R^d)}
 +\omega_{\mathrm{osc}}(t)\\
 &\quad +C_{P,\Omega}
   \|g_h(t)-\hat g_h(t)\|_{L^2(\Omega)} +C_{P,\Omega}
   \|\hat g_h(t)
       -\operatorname{div}\sigma_h(t)\|_{L^2(\Omega)}.
\end{align*}
Taking the weighted $L^2(0,\Tend)$-norm and applying Minkowski's inequality yields
\begin{align*}
 \|R_N\|_{L^2_\beta(0,\Tend;V^*)}
 \le
 \eta_{\mathrm{flux}}(\sigma_h)
 +\eta_{\mathrm{osc}} +\eta_{\mathrm{mean}}
 +\hat\eta_{\mathrm{div}}.
\end{align*}
By \Cref{prop:cdf-jump},
\begin{align*}
 \|\K_{\mu_N-\mu}\tilde u_h\|
 _{L^2_\beta(0,\Tend;V^*)}
 \le\eta_{\mathrm{ker}}^{\mathrm{cdf}}.
\end{align*}
The residual decomposition $R=R_N+\K_{\mu_N-\mu}\tilde u_h$ and the weighted energy estimate in the proof of \Cref{thm:fem-certificate} give \eqref{eq:computable-majorant}.
\end{proof}

\begin{remark}[An explicit domain constant]
\label[remark]{rem:explicit-cp}
If $\Omega$ is contained in an axis-aligned box $\prod_{i=1}^d(a_i,b_i)$, then \eqref{eq:dirichlet-cp} holds with $C_{P,\Omega}=\min_i(b_i-a_i)/\pi$: extend $v$ by zero to the box and apply the one-dimensional Dirichlet Poincar\'e inequality along a coordinate direction.  This also holds for non-convex $\Omega$, and summing over the directions gives the sharper value $[\pi^2\sum_{i=1}^d(b_i-a_i)^{-2}]^{-1/2}$.
\end{remark}

The right-hand side of \eqref{eq:computable-majorant} consists of exact norms and integrals. A verified upper bound additionally requires rigorous upper bounds for every estimator term, including quadrature and algebraic errors, and a certified positive lower bound for $d_{\mu,\beta}$; floating-point quadrature alone does not provide this.

\begin{corollary}[An unweighted certificate beyond the mass-subcritical regime]
\label[corollary]{cor:weight-optimised-certificate}
Retain the problem, measure, reconstruction and flux assumptions of \Cref{prop:computable-majorant}, but do not fix the weight parameter in advance. In particular, $\mu\in\M_+^0([0,\tau])$ has no atom at zero. No mass-subcritical condition $\Lambda_1M<\alpha_0$ is imposed. Define the set of admissible weights by
\begin{align*}
 \mathcal A_\mu:=\{\beta\ge0:\ d_{\mu,\beta}>0\}.
\end{align*}
For each $\beta\in\mathcal A_\mu$, set
\begin{align*}
 \mathcal B_{\beta}^{\mathrm{rec}}
 &:=\frac{\eta_{\mathrm{init}}}{\sqrt{2d_{\mu,\beta}}}
 +\frac{\eta_{\mathrm{flux}}(\beta)+\eta_{\mathrm{osc}}(\beta)
 +\eta_{\mathrm{mean}}(\beta)+\hat\eta_{\mathrm{div}}(\beta)
 +\eta_{\mathrm{ker}}^{\mathrm{rec}}(\beta)}{d_{\mu,\beta}}.
\end{align*}
Then, $\mathcal A_\mu$ is non-empty and
\begin{align}\label{eq:unweighted-optimised}
 \|u_\mu-u_h\|_{L^2(0,\Tend;V)}
 \le\inf_{\beta\in\mathcal A_\mu}
 e^{\beta\Tend}\mathcal B_{\beta}^{\mathrm{rec}}.
\end{align}
For each admissible $\beta$, the quantity $\mathcal B_\beta^{\mathrm{rec}}$ can be evaluated from the data and the reconstruction. If $\supp\mu\subset[s_0,\tau]$ for some $s_0>0$, the explicit sufficient condition
\begin{align}\label{eq:explicit-admissible-beta}
 \Lambda_1M e^{-\beta s_0}<\alpha_0
\end{align}
ensures $\beta\in\mathcal A_\mu$.
\end{corollary}

\begin{proof}
Dominated convergence and $\mu(\{0\})=0$ give $\kappa_\mu(\beta)=\int e^{-\beta s}\,\dd\mu(s)\to0$ as $\beta\to\infty$, proving non-emptiness. For $\beta\in\mathcal A_\mu$, the proof of \Cref{prop:computable-majorant}, combined with \Cref{prop:cancellation-residual}, yields $\|u_\mu-u_h\|_{L^2_\beta(0,\Tend;V)} \le\mathcal B_\beta^{\mathrm{rec}}$. Because $e^{-\beta t}\ge e^{-\beta\Tend}$ on $[0,\Tend]$, \eqref{eq:unweighted-optimised} follows by taking the infimum. The support condition gives $\kappa_\mu(\beta)\le M e^{-\beta s_0}$, which proves \eqref{eq:explicit-admissible-beta}.
\end{proof}

\begin{remark}[Conditioning of the unweighted certificate]
\label[remark]{rem:weight-conditioning}
The factors involving $e^{\beta\Tend}$ and $d_{\mu,\beta}^{-1}$ may be large. The choice of $\beta$ can be optimised a posteriori, since the weighted residual terms depend on the computed reconstruction. Although \eqref{eq:unweighted-optimised} provides a certificate for each admissible delay measure without an atom at zero, the resulting constants are not uniform with respect to the delay locations. To see this, let $\mu=M\delta_{s_0}$ with $\Lambda_1M>\alpha_0$. Then,
\begin{align*}
 d_{\mu,\beta}
 =\alpha_0-\Lambda_1Me^{-\beta s_0},
\end{align*}
and $\beta\in\mathcal A_\mu$ requires
\begin{align*}
 \beta>
 \frac{1}{s_0}
 \log\left(\frac{\Lambda_1M}{\alpha_0}\right).
\end{align*}
Consequently, for every admissible $\beta$,
\begin{align*}
 \frac{e^{\beta\Tend}}{d_{\mu,\beta}}
 >
 \frac{1}{\alpha_0}
 \left(\frac{\Lambda_1M}{\alpha_0}\right)
 ^{\Tend/s_0},
\end{align*}
which diverges as $s_0\downarrow0$. Thus, optimisation over $\beta$ does not remove this loss of uniformity.
\end{remark}

\subsection{Transport budgets and optimal delay atomisation}
\label{subsec:adaptive-delay}
The reliability estimate does not require the delay intervals to have equal length. Let $\{I_j\}_{j=1}^N$ be a measurable partition of $(0,\tau]$, and set $m_j:=\mu(I_j)$. For each $j$ with $m_j>0$, we choose $s_j\in I_j$. We define
\begin{align}\label{eq:adaptive-atom-plan}
 \mu_{\mathcal P}
 &:=\sum_{j:m_j>0}m_j\delta_{s_j},\\
 Q_{\mathcal P}
 &:=\sum_{j:m_j>0}
       \int_{I_j}|s-s_j|\,\dd\mu(s).
\end{align}
The measures $\mu_{\mathcal P}$ and $\mu$ have the same mass $M$, and the coupling that moves $\mu|_{I_j}$ to $s_j$ gives
\begin{align*}
 W_1(\mu,\mu_{\mathcal P})
 \le Q_{\mathcal P}
 \le\sum_{j=1}^N m_j\operatorname{diam}(I_j).
\end{align*}
The cost $Q_{\mathcal P}$ can be evaluated or rigorously bounded above when the corresponding integrals with respect to $\mu$ are available. Under the assumptions of \Cref{thm:fem-certificate}, with $\mu_N=\mu_{\mathcal P}$, the definition \eqref{eq:eta-ker} therefore gives
\begin{align}\label{eq:adaptive-kernel-budget}
 \eta_{\mathrm{ker}}
 \le\Lambda_1\left[
 Q_{\mathcal P}D_h
 +\sqrt{M Q_{\mathcal P}}\|J_h\|_V
 \right]
\end{align}
for any partition, with no mesh-shape-dependent constant.

\begin{proposition}[One-dimensional transport-optimal placement of delay atoms]
\label[proposition]{prop:transport-optimal-atom}
Let $\mu\in\M_+^0([0,\tau])$ have mass $M>0$, and write $Q_\mu(u):=\inf\{s\in[0,\tau]:\mu([0,s])\ge u\}$ for $0<u<M$. Fix positive atomic masses $m_1,\ldots,m_N$ with $\sum_{j=1}^Nm_j=M$, and set $q_0=0$, $q_j=\sum_{i=1}^jm_i$.  Among ordered locations $0<s_1\le\cdots\le s_N\le\tau$, the choice
\begin{align}\label{eq:optimal-atom-quantiles}
 s_j^*:=Q_\mu\left(\frac{q_{j-1}+q_j}{2}\right)
\end{align}
minimises the transport cost, with
\begin{align}\label{eq:optimal-atom-fixed-mass}
 \min_{0<s_1\le\cdots\le s_N\le\tau}
 W_1\left(\mu,\sum_{j=1}^Nm_j\delta_{s_j}\right)
 =\sum_{j=1}^N\int_{q_{j-1}}^{q_j}
       |Q_\mu(u)-s_j^*|\,\dd u.
\end{align}
For any $\beta\ge0$ and $\Tend>0$, the same locations also minimise $\|\Delta_{\nu,\mu}\|_{L^2_\beta(0,\Tend)}$ among such atomic measures $\nu$. If $\dd\mu(s)=(M/\tau)\,\dd s$ on $(0,\tau]$, then, even when both the locations and the masses are free, the unique minimising measure supported on at most $N$ points is
\begin{align}\label{eq:optimal-atom-uniform}
 \mu_N^*=\frac{M}{N}\sum_{j=1}^N
       \delta_{(j-1/2)\tau/N},\quad
 \min_{\substack{\nu\ge0,\ \nu([0,\tau])=M\\
                  |\supp\nu|\le N}}
       W_1(\mu,\nu)=\frac{M\tau}{4N}.
\end{align}
If $\Tend\ge\tau$ and $\beta=0$, this same measure uniquely minimises the cumulative-mass discrepancy, with
\begin{align}\label{eq:optimal-atom-cdf}
 \min_{\substack{\nu\ge0,\ \nu([0,\tau])=M\\
                  |\supp\nu|\le N}}
 \|\Delta_{\nu,\mu}\|_{L^2(0,\Tend)}
 =\frac{M}{N}\sqrt{\frac{\tau}{12}}.
\end{align}
\end{proposition}

\begin{proof}
Monotone coupling gives $W_1(\mu,\nu)=\int_0^M|Q_\mu(u)-Q_\nu(u)|\,\dd u$, where $Q_\nu(u)=s_j$ for $q_{j-1}<u<q_j$, so each interval integral is minimised at a median of $Q_\mu$; \eqref{eq:optimal-atom-quantiles} is one such choice, and it is positive because $\mu(\{0\})=0$. This is classical one-dimensional quantisation~\cite{GrafLuschgy2000,Santambrogio2015}.

For the cumulative-mass discrepancy, put $F_\sigma(t):=\sigma([0,t])$ and $a_j:=(q_{j-1}+q_j)/2$. The defining property of the generalised inverse gives $s_j^*\le t$ if and only if $a_j\le F_\mu(t)$, for $0\le t\le\tau$. Therefore, for $\nu^*:=\sum_{j=1}^Nm_j\delta_{s_j^*}$, the value $F_{\nu^*}(t)$ is a nearest point to $F_\mu(t)$ in $\{q_0,\ldots,q_N\}$ (ties may be resolved either way). For any other ordered placement with the same masses, $F_\nu(t)$ also lies in this set, so $|F_{\nu^*}(t)-F_\mu(t)|\le|F_\nu(t)-F_\mu(t)|$ for $0\le t\le\tau$. Multiplying the squared inequality by $e^{-2\beta t}$ and integrating proves the $L^2_\beta$ minimality, without assuming that $\mu$ has a density.

For the uniform law, $Q_\mu(u)=\tau u/M$. For fixed positive masses $m_1,\ldots,m_N$, the optimal locations give
\begin{align*}
 \min W_1(\mu,\nu)
 =\frac{\tau}{4M}\sum_{j=1}^N m_j^2.
\end{align*}
When $\Tend\ge\tau$ and $\beta=0$, the pointwise nearest-level property gives
\begin{align*}
 \min\|\Delta_{\nu,\mu}\|_{L^2(0,\Tend)}^2
 &=\frac{\tau}{M}
   \sum_{j=1}^N
   2\int_0^{m_j/2}z^2\,\dd z =\frac{\tau}{12M}\sum_{j=1}^N m_j^3.
\end{align*}
By convexity,
\begin{align*}
 \sum_{j=1}^N m_j^2\ge\frac{M^2}{N},
 \quad
 \sum_{j=1}^N m_j^3\ge\frac{M^3}{N^2},
\end{align*}
with equality in either inequality only when $m_j=M/N$ for every $j$. The optimal locations are then uniquely determined by the quantile formula. This proves \eqref{eq:optimal-atom-uniform} and \eqref{eq:optimal-atom-cdf}, including uniqueness. Measures supported on fewer than $N$ points cannot attain these minima.
\end{proof}

\begin{corollary}[Quantile budget with history jumps]
\label[corollary]{cor:quantile-first-order}
Let $\mu\in\M_+^0([0,\tau])$ have mass $M>0$.  For $N\ge1$, we define
\begin{align}\label{eq:equal-mass-quantile-atoms}
 \mu_N:=\frac{M}{N}\sum_{j=1}^N
 \delta_{Q_\mu((j-1/2)M/N)}.
\end{align}
Then, for $0\le t\le\tau$,
\begin{align}\label{eq:quantile-cdf-sup}
 |\Delta_{\mu_N,\mu}(t)|\le\frac{M}{2N},
 \quad W_1(\mu,\mu_N)\le\frac{M\tau}{2N},
\end{align}
and
\begin{align}\label{eq:quantile-cdf-L2}
 \|\Delta_{\mu_N,\mu}\|_{L^2_\beta(0,\Tend)}
 \le\frac{M}{2N}
 \left(\int_0^{\min(\Tend,\tau)}e^{-2\beta t}\,\dd t\right)^{1/2}.
\end{align}
Under the assumptions of \Cref{thm:fem-certificate}, with this choice of $\mu_N$, the delay contribution to the reliability bound satisfies
\begin{align}\label{eq:quantile-certificate-rate}
 \frac{\eta_{\mathrm{ker}}^{\mathrm{cdf}}}{d_{\mu,\beta}}
 \le\frac{\Lambda_1 M}{2N d_{\mu,\beta}}
 \left[\tau D_h+
  \left(\int_0^{\min(\Tend,\tau)}e^{-2\beta t}\,\dd t\right)^{1/2}
  \|J_h\|_V\right].
\end{align}
If, in addition, $\mu$ is uniform, $\Lambda_1M<\alpha_0$, $\beta=0$ and $\Tend\ge\tau$, then, for any \emph{fixed} reconstruction $u_h$, the least possible cdf-based delay contribution among measures supported on at most $N$ atoms equals
\begin{align}\label{eq:optimal-fixed-certificate}
 \min_{\substack{\nu\ge0,\ \nu([0,\tau])=M\\
                  |\supp\nu|\le N}}
 \frac{\eta_{\mathrm{ker}}^{\mathrm{cdf}}(u_h,\nu)}{d_{\mu,0}}
 =\frac{\Lambda_1M}{N d_{\mu,0}}
 \left[\frac{\tau}{4}D_h+
       \sqrt{\frac{\tau}{12}}\|J_h\|_V\right],
\end{align}
attained by the midpoint atoms in \eqref{eq:optimal-atom-uniform}. Here, $\eta_{\mathrm{ker}}^{\mathrm{cdf}}(u_h,\nu)$ means \eqref{eq:cdf-kernel-estimator} evaluated with this fixed reconstruction and the candidate measure $\nu$.
\end{corollary}

\begin{proof}
We write $F_\sigma(t):=\sigma([0,t])$ for $\sigma=\mu,\mu_N$. For $0\le t\le\tau$, the number of quantiles $Q_\mu((j-1/2)M/N)$ not exceeding $t$ equals the number of midpoints $(j-1/2)M/N$ not exceeding $F_\mu(t)$. Consequently, $F_{\mu_N}(t)$ is a nearest multiple of $M/N$ to $F_\mu(t)$, and \eqref{eq:quantile-cdf-sup} follows from the one-dimensional identity $W_1(\mu,\nu)=\int_0^\tau|F_\mu(t)-F_\nu(t)|\,\dd t$. The discrepancy vanishes for $t>\tau$, so integration gives \eqref{eq:quantile-cdf-L2}.  Substitution in \eqref{eq:cdf-kernel-estimator} proves \eqref{eq:quantile-certificate-rate}.  For the uniform law, \Cref{prop:transport-optimal-atom} shows that the midpoint atomic measure minimises both $W_1(\mu,\nu)$ and $\|\Delta_{\nu,\mu}\|_{L^2(0,\Tend)}$ when $\Tend\ge\tau$ and $\beta=0$. Their respective minimum values are
\begin{align*}
 \frac{M\tau}{4N}
 \quad\text{and}\quad
 \frac{M}{N}\sqrt{\frac{\tau}{12}}.
\end{align*}
Because the reconstruction $u_h$ is fixed, $D_h$ and $\|J_h\|_V$ do not depend on the candidate measure $\nu$. The two non-negative terms in $\eta_{\mathrm{ker}}^{\mathrm{cdf}}(u_h,\nu)$ are therefore minimised simultaneously by the midpoint atomic measure. Substitution gives \eqref{eq:optimal-fixed-certificate}.
\end{proof}

\begin{remark}[Delay-law budgets and a posteriori acceptance]
\label[remark]{rem:atom-acceptance}
Fix $\beta\ge0$ with $d_{\mu,\beta}>0$ and a tolerance $\varepsilon_{\mathrm{ker}}>0$. For a general partition, \eqref{eq:adaptive-kernel-budget} shows that the delay-law contribution to the reliability bound does not exceed $\varepsilon_{\mathrm{ker}}$ if
\begin{align*}
 \Lambda_1\left[
 Q_{\mathcal P}D_h
 +\sqrt{M Q_{\mathcal P}}\|J_h\|_V
 \right]
 \le d_{\mu,\beta}\varepsilon_{\mathrm{ker}}.
\end{align*}
For the equal-mass quantile law \eqref{eq:equal-mass-quantile-atoms}, \Cref{cor:quantile-first-order} gives the sufficient condition
\begin{align}\label{eq:quantile-acceptance-test}
 \frac{\Lambda_1 M}{2N d_{\mu,\beta}}
 \left[
 \tau D_h+
 \left(
 \int_0^{\min(\Tend,\tau)}
 e^{-2\beta t}\,\dd t
 \right)^{1/2}
 \|J_h\|_V
 \right]
 \le\varepsilon_{\mathrm{ker}}.
\end{align}
Both conditions involve $D_h$ and $J_h$ and are therefore checked after the reconstruction has been computed; they are acceptance tests, not a priori formulas for $N$. Changing the atoms and recomputing the solution generally changes $D_h$, $J_h$ and the other terms of the majorant. If a test fails, one may refine the atomic approximation and repeat the test with the new reconstruction. For the same reason, the factor $N^{-1}$ in \eqref{eq:quantile-certificate-rate} bounds the delay-law residual on a given reconstruction only. Since $D_h$ and $J_h$ need not remain bounded as the atoms are refined, and since the flux, oscillation, algebraic and initial-data terms of \Cref{prop:computable-majorant} must also be controlled, it does not imply an $N^{-1}$ rate for the total finite element error. For general delay laws, quantile placement is optimal for prescribed atomic masses; optimising the masses is a separate problem~\cite{GrafLuschgy2000}.
\end{remark}

\subsection{A backward Euler realisation}
\label{subsec:backward-euler}
We now state the scheme used in the numerical experiments. Let $0=t_0<t_1<\cdots<t_L=\Tend$, $k_n=t_n-t_{n-1}$, and let $\mu_{\mathcal P}=\sum_{j=1}^Nm_j\delta_{s_j}$ with $m_j\ge0$ and $s_j\in(0,\tau]$. Given $u_h^0\in V_h$, let $u_h$ be the continuous, piecewise affine interpolant of the as-yet-determined states $u_h^n$ for $t\ge0$, with $u_h(0)=u_h^0$. We set $\tilde u_h(t)=\psi(t)$ for $t<0$ and $\tilde u_h(t)=u_h(t)$ for $t\ge0$.  With $f^n:=k_n^{-1}\int_{t_{n-1}}^{t_n}f(t)\,\dd t$, seek $u_h^n\in V_h$ such that
\begin{align}\label{eq:backward-euler-fem}
 &\left(\frac{u_h^n-u_h^{n-1}}{k_n},v_h\right)_H
   +a_0(u_h^n,v_h) \notag \\
   &\quad +\sum_{j=1}^Nm_ja_1\bigl(
   \tilde u_h(t_n-s_j),v_h\bigr)
 =\langle f^n,v_h\rangle_{V^*,V},\quad v_h\in V_h.
\end{align}
When $0<s_j<k_n$, the interpolated delayed state may contain $u_h^n$, so \eqref{eq:backward-euler-fem} is understood as an implicit linear equation.  The delay locations need not be integer multiples of the time step.

\begin{proposition}[Stepwise solvability]
\label[proposition]{prop:be-solvability}
Assume the diffusion model and data hypotheses of \Cref{thm:fem-certificate}, with $\mu_N=\mu_{\mathcal P}$ of mass $M$, and suppose that either $\Lambda_1M<\alpha_0$ or $s_j\ge k_n$ for all $j$ and $n$.  Then, for any conforming finite element space $V_h$ on ${\mathbb T}_h$, any positive time steps $k_n$, and any given $u_h^0\in V_h$, scheme \eqref{eq:backward-euler-fem} has exactly one solution $u_h^n\in V_h$ at each time step. Its piecewise affine reconstruction satisfies \eqref{eq:fem-reconstruction}, and hence admits the reliability bound \eqref{eq:computable-majorant} for every $\beta$ with $d_{\mu,\beta}>0$ and any chosen $L^2$-in-time RT flux reconstruction.
\end{proposition}

\begin{proof}
We proceed step by step. At time $t_n$, all previously computed states and the prescribed history are known.

For each $j$, we set
\begin{align*}
 \theta_{j,n}:=
 \begin{cases}
  1-s_j/k_n,&0<s_j<k_n,\\
  0,&s_j\ge k_n.
 \end{cases}
\end{align*}
If $0<s_j<k_n$, piecewise affine interpolation gives
\begin{align*}
 \tilde u_h(t_n-s_j)
 =\theta_{j,n}u_h^n
 +(1-\theta_{j,n})u_h^{n-1}.
\end{align*}
Otherwise, the delayed value is already known. Thus, the bilinear form for the unknown state $u_h^n$ is
\begin{align*}
 b_n(w,v)
 :=k_n^{-1}(w,v)_H+a_0(w,v) +\sum_{j=1}^Nm_j\theta_{j,n}a_1(w,v).
\end{align*}
Because $m_j\ge0$, $0\le\theta_{j,n}\le1$, and $\sum_{j=1}^Nm_j=M$,
\begin{align*}
 b_n(v,v)
 \ge k_n^{-1}\|v\|_H^2 +(\alpha_0-\Lambda_1M)\|v\|_V^2.
\end{align*}
The mass-subcritical condition makes $b_n$ coercive on $V_h$. If $s_j\ge k_n$ for every $j$, then $\theta_{j,n}=0$, and $b_n(v,v)\ge k_n^{-1}\|v\|_H^2+\alpha_0\|v\|_V^2$ without any condition on $M$. The finite-dimensional Lax--Milgram lemma therefore gives a unique solution at each time step.

The resulting continuous, piecewise affine interpolant belongs to $H^1(0,\Tend;V)$ and takes values in $V_h$. It therefore satisfies \eqref{eq:fem-reconstruction}. The reliability estimate \eqref{eq:computable-majorant} then applies with any $L^2$-in-time RT flux reconstruction, with its equilibration defect included in the majorant.
\end{proof}

\section{A conditional efficiency result}\label{sec:efficiency}
Reliability holds without any shape-regularity restriction, but it says nothing about how close the estimator is to the error. We study efficiency for scalar constant coefficients. The positive result rests on an elliptic equilibration hypothesis, and \Cref{prop:eq-obstruction} shows that this hypothesis cannot hold uniformly for RT fluxes on anisotropic mesh families. Efficiency is measured against a combined error rather than against $u_\mu-u_h$ alone, since cancellation between the delay-law and discretisation errors rules out an unconditional efficiency bound in terms of the latter. The proof uses the cancellation-preserving estimator $\eta_{\rm ker}^{\rm rec}$; the cumulative-mass and transport bounds of \Cref{sec:certified-fem} are intended for choosing the delay atoms.

Assume $B_0=b_0 I$, $B_1=b_1 I$, with $b_0>0$, $b_1>0$, $\Omega$ polyhedral, and a uniform time grid $t_n=nk$. We write $I_n:=(t_{n-1},t_n)$. We retain the non-negative equal-mass, zero-origin-atom hypotheses for $\mu$ and $\mu_N$ from \Cref{thm:fem-certificate}. In this section only, assume the atoms of $\mu_N$ satisfy $s_j=p_jk\in(0,\tau]$, $p_j\in\mathbb N$, and $m_j>0$. The reliability estimate of \Cref{thm:fem-certificate} does not require this alignment; it is used here for the time-residual decomposition and in the counterexample. We further assume $f\in L^2(0,\Tend;L^2(\Omega))$. The finite element states are defined by \eqref{eq:backward-euler-fem}, which is uniquely solvable by \Cref{prop:be-solvability} because $s_j\ge k$; the prescribed history belongs to $H^1(-\tau,0;V)$. We use the conforming piecewise affine space $V_h=\mathbb P^1(\mathbb T_h)\cap H_0^1(\Omega)$; this is needed for the Galerkin orthogonality after projecting the flux onto elementwise constants. We put $f^n=k^{-1}\int_{I_n}f(t)\,\dd t$, $f_k|_{I_n}=f^n$, and
\begin{align*}
 g^n&=(u_h^n-u_h^{n-1})/k-f^n,\\
 F^n&=b_0\nabla u_h^n+b_1\sum_jm_j\nabla\tilde u_h(t_n-s_j),\quad q^n=\Pi_h^0g^n.
\end{align*}
We write $\mathcal E_\beta(w):=(\norm{w}_{L^2_\beta(0,\Tend;V)}^2+\norm{\partial_tw}_{L^2_\beta(0,\Tend;V^*)}^2)^{1/2}$ and
$\mathcal E_{\rm comb}:=\mathcal E_\beta(u_\mu-u_h)+\mathcal E_\beta(u_{\mu_N}-u_h)$.
In this scalar setting we take $\alpha_0=\Lambda_0=b_0$ and $\Lambda_1=b_1$. We require $d_\nu=b_0-b_1\kappa_\nu(\beta)>0$ for $\nu=\mu,\mu_N$, and put $L_\nu=b_0+b_1\kappa_\nu(\beta)$.

Let $C_P$ be any Dirichlet Poincar\'e constant, so that $\norm{v}_H\le C_P\norm{v}_V$ for $v\in V$. For $q\in L^2(\Omega)$, we write
\begin{align*}
 \mathcal Q_h(q):=
 \{\sigma_h\in\mathbb{RT}^0(\mathbb T_h):
   \operatorname{div}\sigma_h=\Pi_h^0q\}.
\end{align*}

\begin{lemma}[Minimal RT fluxes]\label[lemma]{lem:crrt-optimality}
For $q\in L^2(\Omega)$, the set $\mathcal Q_h(q)$ contains a unique element $\sigma_h^M(q)$ of minimal $L^2$-norm. For any $z\in V$ and $\sigma_h\in\mathcal Q_h(q)$,
\begin{align}\label{eq:crrt-min-identity}
 \norm{\nabla z-\sigma_h}_{L^2(\Omega)^d}^2
 =\norm{\nabla z-\sigma_h^M(q)}_{L^2(\Omega)^d}^2
 +\norm{\sigma_h-\sigma_h^M(q)}_{L^2(\Omega)^d}^2.
\end{align}
Thus, $\sigma_h^M(q)$ minimises $\norm{\nabla z-\sigma_h}_{L^2(\Omega)^d}$ over $\mathcal Q_h(q)$ for every $z\in V$, without any condition on the mesh geometry. It is the lowest-order RT mixed flux for the homogeneous Dirichlet problem with divergence datum $\Pi_h^0q$, and can be computed by Marini's Crouzeix--Raviart post-processing~\cite{Marini1985,ArnoldBrezzi1985}.
\end{lemma}

\begin{proof}
By \Cref{subsec:fem-reliability}, $\mathcal Q_h(q)$ is a non-empty affine subspace of the finite-dimensional space $\mathbb{RT}^0(\mathbb T_h)$. Therefore, $\sigma_h^M(q)$ exists, is unique, and is orthogonal to every divergence-free $\delta_h\in\mathbb{RT}^0(\mathbb T_h)$. For $z\in V$, integration by parts gives $(\nabla z,\delta_h)=-(z,\operatorname{div}\delta_h)=0$. Therefore, $\nabla z-\sigma_h^M(q)$ is orthogonal to $\sigma_h-\sigma_h^M(q)$, and the square expansion gives \eqref{eq:crrt-min-identity}.
\end{proof}

An anisotropic CR--RT setting for semilinear elliptic problems is discussed in~\cite{IshizakaCRRT2026}, while the related scalar-matrix derivation is given in~\cite{IshizakaMariniL2_2026}.

In the scalar setting, $F_N(t)$ and $F^n$ are gradients of functions in $V$. Therefore, for $t\in I_n$, the flux $\sigma_h^M(q^n)$ minimises both $\norm{F_N(t)-\sigma_h}_{L^2(\Omega)^d}$ and $\norm{F^n-\sigma_h}_{L^2(\Omega)^d}$ over $\mathcal Q_h(g^n)$. This does not identify the minimiser for the time-dependent datum $\Pi_h^0g(t)$ with that for $q^n$.

We say that the spatial efficiency condition (EQ) holds for the step data if there exists a constant $C_{\rm eq}\ge1$, independent of the time-step index, such that
\begin{align}
 \inf_{\sigma_h\in\mathcal Q_h(q)}
 \norm{F_h-\sigma_h}_{L^2(\Omega)^d}
 &\le C_{\rm eq}
 \left(\norm{\ell}_{V^*}+\mathrm{osc}_h(q)\right),
 \tag{EQ}\label{eq:elliptic-EQ}\\
 \mathrm{osc}_h(q)^2
 &:=\sum_T\frac{h_T^2}{\pi^2}
 \norm{q-\Pi_h^0q}_{L^2(T)}^2,
 \notag
\end{align}
for the pairs $(q,F_h)=(g^n,\Pi_h^0F^n)$, where
\begin{align*}
 \ell(v):=-(q,v)-(F_h,\nabla v),
 \quad v\in V.
\end{align*}
These pairs satisfy $\ell(v_h)=0$ for any $v_h\in V_h$. The constant $C_{\rm eq}$ may depend on the spatial mesh; by \Cref{prop:eq-obstruction}, it cannot be bounded uniformly over general anisotropic mesh families. The identity \eqref{eq:crrt-min-identity} alone does not imply \eqref{eq:elliptic-EQ}.

\begin{remark}[Scope of (EQ)]\label[remark]{rem:eq-shapereg}
Suppose that $F_h=\nabla z_h$ with $z_h\in V_h$ and that $\ell$ vanishes on $V_h$. If $z\in V$ solves $(\nabla z,\nabla v)=-(q,v)$ for $v\in V$, then $\ell(v)=(\nabla(z-z_h),\nabla v)$, so $\norm{\ell}_{V^*}=\norm{\nabla(z-z_h)}_{L^2(\Omega)^d}$. For the two-dimensional Poisson problem on shape-regular triangulations, vertex-patch flux equilibration and local efficiency are established in \cite{BraessPillweinSchoeberl2009,ErnVohralik2015}. The precise form of (EQ) additionally fixes the global $\mathbb{RT}^0$ space, the projected divergence and the oscillation normalisation, so we do not invoke these results as a proof of (EQ); \Cref{thm:conditional-efficiency} applies whenever a valid $C_{\rm eq}$ is available for the mesh and data. In the present scheme, if all history samples used at the time nodes belong to $V_h$, then $F^n$ is the gradient of a function in $V_h$.
\end{remark}

We define the stepwise flux $\sigma_h(t)=\sigma_h^M(q^n)$ on $I_n$. Its divergence is $\Pi_h^0g^n$ rather than $\Pi_h^0g(t)$, so it need not satisfy \eqref{eq:equilibrated-flux}; in this section, $\eta_{\rm flux}$ denotes $\eta_{\rm flux}(\sigma_h)$ from \eqref{eq:eta-flux-sigma} for this flux. We define
\begin{align*}
 \eta_{f,k}&:=C_P\norm{f-f_k}_{L^2_\beta(0,\Tend;L^2(\Omega))},\\
 \eta_{\mathrm{osc},k}^2&:=\sum_n\int_{I_n}e^{-2\beta t}\sum_T\frac{h_T^2}{\pi^2}\norm{g^n-\Pi_h^0g^n}_{L^2(T)}^2\,\dd t,\\
 \eta_{F,k}^2&:=\sum_n\int_{I_n}e^{-2\beta t}\norm{F^n-\Pi_h^0F^n}_{L^2(\Omega)^d}^2\,\dd t,\\
 \eta_{\psi,k}&:=b_1M\norm{\psi-\psi_k}_{L^2(-p_{\max}k,0;V)},\\
 \eta_{J,k}^2&:=b_1^2\norm{J_h}_V^2\sum_n\int_{I_n}e^{-2\beta t}
 \Bigl(\sum_{j:p_j=n}m_j\Bigr)^2\,\dd t,\\
 \mathrm{osc}_{\rm tot}&:=\eta_{\mathrm{osc},k}+\eta_{f,k}+\eta_{F,k}+\eta_{\psi,k}+\eta_{J,k},\\
 \Sigma_\beta&:=\eta_{\rm flux}+\eta_{\rm ker}^{\rm rec}(\beta)+\eta_{\rm init}.
\end{align*}
Here, $p_{\max}:=\max_{1\le j\le N}p_j$, and $\psi_k$ denotes the continuous, piecewise affine interpolant of the prescribed history at the negative time nodes $-p_{\max}k,\ldots,-k,0$. In particular, $\psi_k(0)=\psi(0)$. The possible jump $J_h=u_h(0)-\psi(0)$ is accounted for separately by $\eta_{J,k}$.

For sufficiently smooth data, $\eta_{\mathrm{osc},k}$ and $\eta_{\psi,k}$ can decay under refinement, and $\eta_{J,k}\le b_1M\sqrt k\,\norm{J_h}_V$. The term $\eta_{f,k}$ measures temporal approximation of the forcing. The term $\eta_{F,k}$ measures the elementwise constant projection error of the flux at the time nodes. It vanishes if all relevant history samples belong to $V_h$. For general history data, its decay under spatial refinement depends on their regularity; no first-order rate is asserted under the present assumptions.

\begin{theorem}[Two-sided bounds for the combined error]\label[theorem]{thm:conditional-efficiency}
Under the assumptions above, let $A_\nu=(1+L_\nu)/\sqrt{2d_\nu}$, $B_\nu=1+(1+L_\nu)/d_\nu$, and $C_\nu^+=(1+L_\nu^2)^{1/2}$. Then,
\begin{align}\label{eq:graph-reliability}
\mathcal E_{\rm comb} \leq (B_\mu+B_{\mu_N})(\eta_{\rm flux}+\eta_{\mathrm{osc},k}+\eta_{f,k})
 +B_\mu\eta_{\rm ker}^{\rm rec}+(A_\mu+A_{\mu_N})\eta_{\rm init}.
\end{align}
If, in addition, \eqref{eq:elliptic-EQ} holds, then
\begin{align}\label{eq:graph-efficiency}
 \Sigma_\beta\le C_{\rm eff} \left(\mathcal E_{\rm comb}+\mathrm{osc}_{\rm tot}\right),
 \quad
 C_{\rm eff}:=\max\bigl\{C_\mu^++C_{\rm tr},\,(1+C_F)C_{\mu_N}^+\bigr\},
\end{align}
where $C_{\rm tr}:=e^{\beta\Tend}(1+C_P^2/\Tend)^{1/2}$, $C_F:=e^{\beta k}\sqrt{2/\lambda_*}\,C_{\rm eq}$, and $\lambda_*:=(4-\sqrt{13})/6$.
\end{theorem}

\begin{proof}
\emph{Step 1: residual bounds and reliability.} For $\nu\in\{\mu,\mu_N\}$, we set $w_\nu:=u_\nu-u_h$ and
\begin{align*}
 R_\nu:=f-\partial_tu_h-A_0u_h-\mathcal K_\nu\tilde u_h.
\end{align*}
The error has zero history and satisfies $\partial_tw_\nu+A_0w_\nu+\mathcal K_\nu\tilde w_\nu=R_\nu$. The weighted energy estimate of \Cref{thm:fem-certificate} and the bound $\|A_0w_\nu+\mathcal K_\nu\tilde w_\nu\|_{L^2_\beta(0,\Tend;V^*)} \le L_\nu\|w_\nu\|_{L^2_\beta(0,\Tend;V)}$ give
\begin{align*}
 \mathcal E_\beta(w_\nu)
 \le B_\nu\norm{R_\nu}_{L^2_\beta(0,\Tend;V^*)}
       +A_\nu\norm{w_\nu(0)}_H,\quad
 \norm{R_\nu}_{L^2_\beta(0,\Tend;V^*)}
 \le C_\nu^+\mathcal E_\beta(w_\nu),
\end{align*}
where $w_\nu(0)=u_0-u_h(0)$ for both delay laws. For $t\in I_n$, the backward Euler reconstruction satisfies $g(t)=g^n+f_k(t)-f(t)$. Because $\operatorname{div}\sigma_h=q^n$ on $I_n$, integration by parts yields, for $v\in V$,
\begin{align*}
 \pair{R_{\mu_N}(t)}{v}
 &=-(g^n-q^n,v)_{L^2(\Omega)}
      +(f(t)-f_k(t),v)_{L^2(\Omega)}\\
 &\quad -(F_N(t)-\sigma_h(t),\nabla v)_{L^2(\Omega)^d}.
\end{align*}
The elementwise Poincar\'e inequality, the global Dirichlet Poincar\'e inequality and Minkowski's inequality give $\norm{R_{\mu_N}}_{L^2_\beta(0,\Tend;V^*)}\le\eta_{\rm flux}+\eta_{\mathrm{osc},k}+\eta_{f,k}$. By \eqref{eq:cancellation-identity}, $R_\mu-R_{\mu_N}=A_1\mathcal C_h$, and therefore $\norm{R_\mu}_{L^2_\beta(0,\Tend;V^*)}\le\norm{R_{\mu_N}}_{L^2_\beta(0,\Tend;V^*)}+\eta_{\rm ker}^{\rm rec}$. These inequalities prove \eqref{eq:graph-reliability}.

\emph{Step 2: initial value.} For $w\in\W(0,\Tend)$, integrating $\frac{\mathrm d}{\mathrm dt}\bigl(\varphi\norm{w}_H^2\bigr)$ with $\varphi(t)=(\Tend-t)/\Tend$ over $(0,\Tend)$ gives
$\norm{w(0)}_H^2\le\Tend^{-1}\norm{w}_{L^2(0,\Tend;H)}^2+2\norm{\partial_tw}_{L^2(0,\Tend;V^*)}\norm{w}_{L^2(0,\Tend;V)}$. The Poincar\'e and Young inequalities and $e^{-2\beta t}\ge e^{-2\beta\Tend}$ on $(0,\Tend)$ yield $\norm{w(0)}_H\le C_{\rm tr}\mathcal E_\beta(w)$; in particular, $\eta_{\rm init}\le C_{\rm tr}\mathcal E_\beta(w_\mu)$.

\emph{Step 3: time decomposition.} For $t\in I_n$, we set $\theta:=(t-t_{n-1})/k$ and $M_n:=\sum_{j:p_j=n}m_j$. For every integer $m$ occurring below, we set $t_m:=mk$ and define
\begin{align*}
 \delta^m:=
 \begin{cases}
  u_h^m-u_h^{m-1},&m\ge1,\\
  \psi(t_m)-\psi(t_{m-1}),&m\le0,
 \end{cases}
 \quad
 w^n:=b_0\delta^n+b_1\sum_{j=1}^Nm_j\delta^{n-p_j}.
\end{align*}
In particular, $\delta^0=\psi(0)-\psi(-k)$, and the jump $J_h$ is kept separate. On an interval for which $t-s_j<0$, we write the history as $\psi_k(t-s_j)$ plus $(\psi-\psi_k)(t-s_j)$. Piecewise affine interpolation then gives
\begin{align*}
 F_N(t)-F^n
 = -(1-\theta)\nabla w^n-b_1M_n\nabla J_h +b_1\sum_{j:p_j\ge n}m_j
       \nabla(\psi-\psi_k)(t-s_j).
\end{align*}
The term involving $J_h$ occurs precisely for $p_j=n$, when the delayed argument crosses zero during $I_n$. Let $\mathcal R_V:V\to V^*$ be the gradient Riesz map, $\pair{\mathcal R_Vv}{z}=(\nabla v,\nabla z)_{L^2(\Omega)^d}$, and define
\begin{align*}
 \ell^n(v)&:=-(g^n,v)_{L^2(\Omega)}
           -(F^n,\nabla v)_{L^2(\Omega)^d},\\
 r_\psi(t)&:=-b_1\sum_{j:p_j\ge n}m_j
            \mathcal R_V[(\psi-\psi_k)(t-s_j)],\\
 R_{\mu_N}^{\flat}(t)&:=R_{\mu_N}(t)-(f(t)-f_k(t)).
\end{align*}
The flux identity gives the exact decomposition in $V^*$
\begin{align*}
 R_{\mu_N}^{\flat}(t)-r_\psi(t)
 =a_n+(1-\theta)b_n,\quad
 a_n:=\ell^n+b_1M_n\mathcal R_VJ_h,\quad
 b_n:=\mathcal R_Vw^n.
\end{align*}
By Minkowski's inequality, $\sum_jm_j=M$ and $\beta\ge0$, $\norm{r_\psi}_{L^2_\beta(0,\Tend;V^*)}\le\eta_{\psi,k}$.

\emph{Step 4: spatial equilibration.} The backward Euler equations give $\ell^n(v_h)=0$ for any $v_h\in V_h$. We put $D^n:=F^n-\Pi_h^0F^n$ and
\begin{align*}
 \bar\ell^n(v):=-(g^n,v)_{L^2(\Omega)}
     -(\Pi_h^0F^n,\nabla v)_{L^2(\Omega)^d}.
\end{align*}
Because $\nabla v_h$ is elementwise constant, $(D^n,\nabla v_h)_{L^2(\Omega)^d}=0$; hence, $\bar\ell^n(v_h)=0$ and $\norm{\bar\ell^n}_{V^*}\le\norm{\ell^n}_{V^*}+\norm{D^n}_{L^2(\Omega)^d}$. By \Cref{lem:crrt-optimality}, $\sigma_h^M(q^n)$ minimises $\norm{F^n-\sigma_h}_{L^2(\Omega)^d}$ over $\mathcal Q_h(g^n)$. The same flux minimises the distance from $\Pi_h^0F^n$: differences of admissible RT fluxes are divergence-free and thus elementwise constant, so they are orthogonal to $D^n$. Applying (EQ) to $(g^n,\Pi_h^0F^n)$ and using the triangle inequality gives
\begin{align*}
 \norm{F^n-\sigma_h^M(q^n)}_{L^2(\Omega)^d}
 \le(1+C_{\rm eq})\norm{D^n}_{L^2(\Omega)^d} +C_{\rm eq}\bigl(\norm{\ell^n}_{V^*}
                  +\mathrm{osc}_h(g^n)\bigr).
\end{align*}

\emph{Step 5: inverse estimate in time.} The Gram matrix of $1$ and $1-\theta$ on $(0,1)$ is
\begin{align*}
 G=\begin{pmatrix}1&1/2\\1/2&1/3\end{pmatrix},
 \quad \lambda_{\min}(G)=\lambda_*.
\end{align*}
Consequently, $\int_0^1\norm{a+(1-\theta)b}_{V^*}^2\,\dd\theta\ge\lambda_*(\norm a_{V^*}^2+\norm b_{V^*}^2)$ for any $a,b\in V^*$. Because $e^{-2\beta t_n}\le e^{-2\beta t} \le e^{2\beta k}e^{-2\beta t_n}$ on $I_n$, Step~3, the Cauchy--Schwarz inequality and $C_{\rm eq}\ge1$ imply
\begin{align*}
 &C_{\rm eq}\left(\sum_n\int_{I_n}e^{-2\beta t}
                   \norm{a_n}_{V^*}^2\,\dd t\right)^{1/2} + \left(\sum_n\int_{I_n}e^{-2\beta t}(1-\theta)^2
                   \norm{b_n}_{V^*}^2\,\dd t\right)^{1/2}\\
 &\quad \le C_F\norm{R_{\mu_N}^{\flat}-r_\psi}_{L^2_\beta(0,\Tend;V^*)}
 \le C_F\bigl(\norm{R_{\mu_N}}_{L^2_\beta(0,\Tend;V^*)}
      +\eta_{f,k}+\eta_{\psi,k}\bigr).
\end{align*}

\emph{Step 6: conclusion.} We have $\norm{\ell^n}_{V^*}\le\norm{a_n}_{V^*}+b_1M_n\norm{J_h}_V$ and $\norm{b_n}_{V^*}=\norm{w^n}_V$. The bounds for $F^n-\sigma_h^M(q^n)$ and $F_N(t)-F^n$, Step~5 and Minkowski's inequality yield
\begin{align*}
 \eta_{\rm flux} &\le
 C_F\norm{R_{\mu_N}}_{L^2_\beta(0,\Tend;V^*)}
 +C_F\eta_{f,k}+(1+C_F)\eta_{\psi,k}\\
 &\quad +(1+C_{\rm eq})(\eta_{F,k}+\eta_{J,k})
  +C_{\rm eq}\eta_{\mathrm{osc},k}.
\end{align*}
In the scalar setting, $A_1=b_1\mathcal R_V$, so \eqref{eq:eta-rec-exact} gives
\begin{align*}
 \eta_{\rm ker}^{\rm rec}
 =\norm{R_\mu-R_{\mu_N}}_{L^2_\beta(0,\Tend;V^*)}
 \le\norm{R_\mu}_{L^2_\beta(0,\Tend;V^*)}
    +\norm{R_{\mu_N}}_{L^2_\beta(0,\Tend;V^*)}.
\end{align*}
Combining these inequalities with Steps~1 and~2, and using $C_{\rm eq}\le C_F$, we obtain
\begin{align*}
 \Sigma_\beta \le
 (C_\mu^++C_{\rm tr})\mathcal E_\beta(w_\mu)
 +(1+C_F)C_{\mu_N}^+\mathcal E_\beta(w_{\mu_N})
 +(1+C_F)\,\mathrm{osc}_{\rm tot}.
\end{align*}
Since $C_{\mu_N}^+\ge1$, this proves \eqref{eq:graph-efficiency}.
\end{proof}

\begin{proposition}[An aspect-ratio obstruction]\label[proposition]{prop:eq-obstruction}
Let $\Omega=(0,1)^2$. For each even integer $p\ge2$, partition $\Omega$ into $p\times p^2$ congruent rectangles of widths $h_x=p^{-1}$ and heights $h_y=p^{-2}$, and split every rectangle along the same diagonal into two right triangles. Let $V_h$ be the conforming homogeneous Dirichlet $\mathbb P^1$ space on this mesh. We number the horizontal strips $S_j:=(0,1)\times(jh_y,(j+1)h_y)$, $j=0,\ldots,p^2-1$, from bottom to top, and set
\begin{align*}
 r_p:=(-1)^j\ \text{on }S_j,\quad q_p:=-r_p,\quad F_p:=0,\quad
 \ell_p(v):=-(q_p,v)-(F_p,\nabla v)=(r_p,v).
\end{align*}
\begin{enumerate}
\renewcommand{\labelenumi}{(\roman{enumi})}
\item $\ell_p(v_h)=0$ for any $v_h\in V_h$, $\mathrm{osc}_h(q_p)=0$, and
\begin{align}\label{eq:EQ-counterexample-residual-bound}
 \frac{h_y}{\pi\sqrt3}\le\norm{\ell_p}_{V^*}\le\frac{h_y}{\sqrt3}.
\end{align}
\item Every $\sigma_h\in\mathcal Q_h(q_p)$ satisfies $\norm{\sigma_h}_{L^2(\Omega)^2}\ge h_x/\sqrt{72}$. Consequently,
\begin{align}\label{eq:EQ-counterexample}
\frac{\displaystyle\inf_{\sigma_h\in\mathcal Q_h(q_p)}
\|F_p-\sigma_h\|_{L^2(\Omega)^2}}
{\|\ell_p\|_{V^*}+\mathrm{osc}_h(q_p)}
\ge \frac{h_x}{\sqrt{24}\,h_y}=\frac{p}{\sqrt{24}}.
\end{align}
\item In the setting of \Cref{thm:conditional-efficiency}, take $u_0=u_h^0=0$, $\psi=0$ and $f(x,t)=r_p(x)$ for $0<t<\Tend$. Then, the backward Euler solution vanishes, $\eta_{\rm ker}^{\rm rec}=\eta_{\rm init}=\mathrm{osc}_{\rm tot}=0$, and, with $\omega_\beta:=(\int_0^{\Tend}e^{-2\beta t}\,\dd t)^{1/2}$,
\begin{align}\label{eq:nonuniform-estimator}
 \Sigma_\beta \ge \omega_\beta\frac{h_x}{\sqrt{72}}, \quad
 \mathcal E_{\rm comb}\le(B_\mu+B_{\mu_N})\omega_\beta\frac{h_y}{\sqrt3},\quad
 \frac{\Sigma_\beta}{\mathcal E_{\rm comb}+\mathrm{osc}_{\rm tot}}
  \ge\frac{p}{(B_\mu+B_{\mu_N})\sqrt{24}}.
\end{align}
\end{enumerate}
Every triangle has maximum angle $\pi/2$. Therefore, neither the maximum-angle condition nor this structured anisotropic refinement ensures an aspect-ratio-uniform constant in \eqref{eq:elliptic-EQ}. Since $B_\mu$ and $B_{\mu_N}$ depend only on $b_0$, $b_1$, $\beta$ and the delay laws, the constant in \eqref{eq:graph-efficiency} grows at least linearly with the aspect ratio $h_x/h_y=p$ for fixed $k$, $\beta$, $\mu$ and $\mu_N$, whereas the reliability bound \eqref{eq:graph-reliability} remains valid.
\end{proposition}

\begin{proof}
\emph{(i)} Because $q_p$ is constant on each triangle, $\Pi_h^0q_p=q_p$, and hence $\mathrm{osc}_h(q_p)=0$. Let $\varphi_h$ be the nodal basis function associated with an interior vertex. Its support meets the horizontal strip immediately above the vertex in three triangles, and the strip immediately below in three triangles. Each triangle has area $h_xh_y/2$, and the integral of $\varphi_h$ over each of them is one-third of its area. Consequently, $\varphi_h$ has integral $h_xh_y/2$ over each of the two strips, on which $r_p$ takes opposite values. Therefore, $(r_p,\varphi_h)_{L^2(\Omega)}=0$, and since the interior nodal basis functions span $V_h$, $\ell_p$ vanishes on $V_h$.

We set $R_p(y):=\int_0^y r_p(x,s)\,\dd s$, which is independent of $x$. On successive strips, $R_p$ increases linearly from $0$ to $h_y$ and decreases linearly from $h_y$ to $0$. Direct integration gives
\begin{align*}
 \norm{R_p}_{L^2(0,1)}^2=\frac{h_y^2}{3},\quad
 \int_0^1R_p(y)\,\dd y=\frac{h_y}{2}.
\end{align*}
For any $v\in H_0^1(\Omega)$, integration by parts in the $y$-direction yields
\begin{align*}
 \ell_p(v)=-\int_\Omega R_p(y)\,\partial_yv(x,y)\,\dd x\,\dd y,
\end{align*}
and the Cauchy--Schwarz inequality gives the upper bound in \eqref{eq:EQ-counterexample-residual-bound}. For the lower bound, let $\chi_p'=h_y/2-R_p$ with $\chi_p(0)=0$. Then, $\chi_p(1)=0$ and $\norm{\chi_p'}_{L^2(0,1)}=h_y/\sqrt{12}$. For $v(x,y):=\sin(\pi x)\chi_p(y)\in H_0^1(\Omega)$, integration by parts and $\int_0^1\chi_p'\,\dd y=0$ give
\begin{align*}
 \ell_p(v)=-\frac{2}{\pi}\int_0^1R_p\chi_p'\,\dd y=\frac{2}{\pi}\norm{\chi_p'}_{L^2(0,1)}^2,
\end{align*}
whereas the one-dimensional Dirichlet Poincar\'e inequality gives
$\norm{\nabla v}_{L^2(\Omega)^2}^2=\frac12(\pi^2\norm{\chi_p}_{L^2(0,1)}^2+\norm{\chi_p'}_{L^2(0,1)}^2)\le\norm{\chi_p'}_{L^2(0,1)}^2$.
Therefore, $\norm{\ell_p}_{V^*}\ge\ell_p(v)/\norm{v}_V\ge(2/\pi)\norm{\chi_p'}_{L^2(0,1)}=h_y/(\pi\sqrt3)$.

\emph{(ii)} The set $\mathcal Q_h(q_p)$ is non-empty by \Cref{lem:crrt-optimality}. On each triangle $T$ with barycentre $x_T$, the local RT representation and the prescribed divergence give $\sigma_h|_T=a_T+\frac{q_p|_T}{2}(x-x_T)$ with $a_T\in\R^2$. Because $\int_T(x-x_T)\,\dd x=0$, the two terms are orthogonal in $L^2(T)^2$. Every triangle is a right triangle with perpendicular side lengths $h_x$ and $h_y$, so a direct calculation of its second moment gives $\int_T|x-x_T|^2\,\dd x=\frac{|T|_d}{18}(h_x^2+h_y^2)$. Because $|q_p|=1$,
\begin{align*}
 \norm{\sigma_h}_{L^2(T)^2}^2
 \ge\frac14\int_T|x-x_T|^2\,\dd x
 =\frac{|T|_d}{72}(h_x^2+h_y^2).
\end{align*}
Summing over all triangles and using $|\Omega|=1$ gives $\norm{\sigma_h}_{L^2(\Omega)^2}\ge h_x/\sqrt{72}$. Combined with \eqref{eq:EQ-counterexample-residual-bound} and $\mathrm{osc}_h(q_p)=0$, this proves \eqref{eq:EQ-counterexample}, since $(h_x/\sqrt{72})/(h_y/\sqrt3)=h_x/(\sqrt{24}\,h_y)$.

\emph{(iii)} Because the atoms are aligned with the time grid and $p_j\ge1$, each delayed value at $t_n$ depends only on the prescribed history or previously computed states. Suppose that $u_h^m=0$ for $0\le m<n$. As $\psi=0$, the step-$n$ equation \eqref{eq:backward-euler-fem} reduces to
\begin{align*}
 k^{-1}(u_h^n,v_h)_H
 +b_0(\nabla u_h^n,\nabla v_h)_{L^2(\Omega)^2}
 =(r_p,v_h)_{L^2(\Omega)}=0,
 \quad v_h\in V_h,
\end{align*}
by part~(i). Because $b_0>0$, coercivity gives $u_h^n=0$. By induction, $u_h=0$ on $[0,\Tend]$, and $J_h=D_h=0$. The discrete data satisfy $g^n=-r_p=q_p=\Pi_h^0g^n$ and $F^n=0$. Because $f$ is constant in time and the history vanishes, all contributions to $\mathrm{osc}_{\rm tot}$ are zero, and $\eta_{\rm ker}^{\rm rec}=\eta_{\rm init}=0$. Therefore, $\Sigma_\beta=\eta_{\rm flux}$. Since $F_N(t)=0$ and $\sigma_h(t)=\sigma_h^M(q_p)$ on every $I_n$, part~(ii) gives
\begin{align*}
 \Sigma_\beta=\omega_\beta\norm{\sigma_h^M(q_p)}_{L^2(\Omega)^2}\ge\omega_\beta\frac{h_x}{\sqrt{72}}.
\end{align*}
For $\nu\in\{\mu,\mu_N\}$, the error $w_\nu=u_\nu$ has vanishing initial value and history, and $R_\nu(t)=r_p$ for $0<t<\Tend$. Part~(i) gives $\norm{R_\nu}_{L^2_\beta(0,\Tend;V^*)}=\omega_\beta\norm{r_p}_{V^*}\le\omega_\beta h_y/\sqrt3$, and the graph estimate in Step~1 of the proof of \Cref{thm:conditional-efficiency} yields $\mathcal E_\beta(w_\nu)\le B_\nu\omega_\beta h_y/\sqrt3$. Summing over $\nu$ gives the bound for $\mathcal E_{\rm comb}$, and the ratio in \eqref{eq:nonuniform-estimator} follows from $h_x/h_y=p$.
\end{proof}

\begin{remark}[Brezzi--Douglas--Marini fluxes and related work]
\label[remark]{rem:bdm-comparison}
The obstruction arises from the elementwise term $(q_T/d)(x-x_T)$ imposed by the local structure of $\mathbb{RT}^0(T)$; it persists even for the least-squares optimal RT flux. The lowest-order Brezzi--Douglas--Marini space $\mathrm{BDM}_1(\mathbb T_h)$ consists of $H(\operatorname{div};\Omega)$-conforming piecewise affine fields of the local form $\boldsymbol a_T+B_T(x-x_T)$ with $B_T\in\R^{2\times2}$. For the strip data of \Cref{prop:eq-obstruction}, the continuous piecewise affine field $\sigma_h^{\mathrm{BDM}}:=(0,-R_p(y))^\top$ belongs to $\mathrm{BDM}_1(\mathbb T_h)$, satisfies $\operatorname{div}\sigma_h^{\mathrm{BDM}}=q_p$ and, by part~(i) of \Cref{prop:eq-obstruction},
\begin{align*}
 \norm{\sigma_h^{\mathrm{BDM}}}_{L^2(\Omega)^2}
 =\frac{h_y}{\sqrt3}\le\pi\norm{\ell_p}_{V^*}.
\end{align*}
The proof of \Cref{thm:fem-certificate} uses only $\sigma_h\in H(\operatorname{div};\Omega)$ with $\operatorname{div}\sigma_h=g_h$, and therefore applies verbatim to such fluxes. This comparison concerns only the constructed mesh and data family; it does not establish an aspect-ratio-uniform efficiency estimate for BDM1 fluxes on general anisotropic meshes.

Grosman~\cite{Grosman2006} showed that the lower bound of the standard equilibrated residual method contains a factor that can be unbounded on anisotropic meshes, and proposed a modified estimator. Kopteva~\cite{Kopteva2017Flux} proposed an anisotropic flux equilibration using corrections across short edges. \Cref{prop:eq-obstruction} gives an explicit instance for a fixed flux space and divergence constraint: every global RT flux with the prescribed divergence, including the least-squares minimiser, satisfies the lower bound \eqref{eq:EQ-counterexample}, although the data oscillation vanishes, all triangles satisfy the maximum-angle condition, and the data arise from the backward Euler discretisation. Neither result rules out other anisotropic flux reconstructions with uniform efficiency.
\end{remark}

%%From here

\section{Numerical experiments}\label{sec:numerics}
We report double-precision quadrature evaluations of the complete majorant under refinement. The first subsection uses $B^{(Q)}$ for \eqref{eq:computable-majorant}, and the grid-aligned tests use its complete stepwise counterpart $B_{\rm step}^{(Q)}$ from \eqref{eq:Bstep}. These values are not validated numerical enclosures: quadrature and rounding errors are not rigorously bounded. We do not compute $\mathcal E_{\rm comb}$, so \eqref{eq:graph-efficiency} is not tested numerically.

\subsection{Subcritical delay law and residual separation}
\begin{figure}[htbp]
\centering
\includegraphics[width=.88\linewidth]{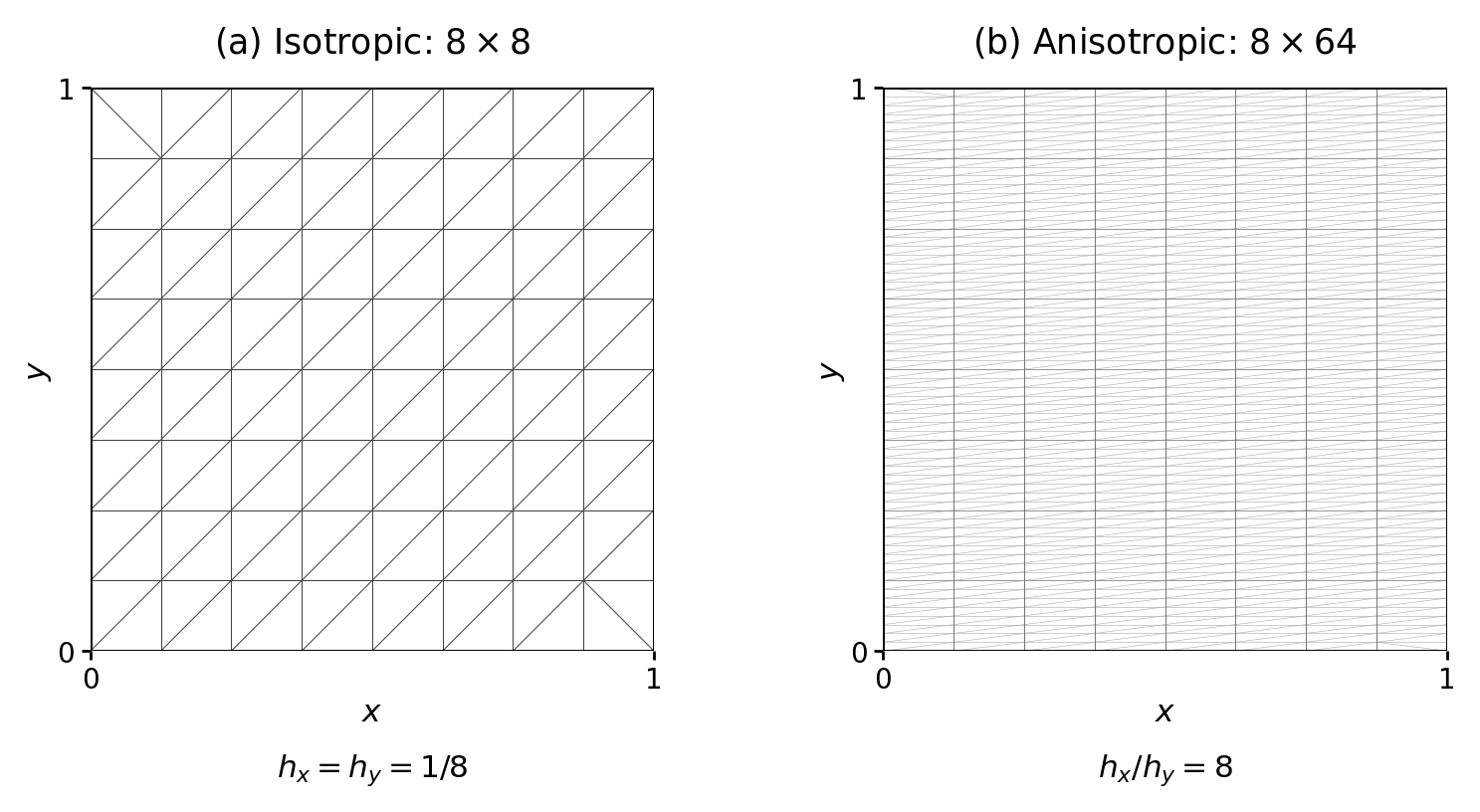}
\caption{Corner-corrected triangulations: (a) the isotropic $8\times8$ mesh of \Cref{tab:manufactured-errors}; (b) the anisotropic $8\times64$ mesh ($h_x/h_y=8$) from the jointly refined family $n_y=n_x^2$.}
\label{fig:corner-mesh}
\end{figure}

Let $\Omega=(0,1)^2$, $\Tend=0.8$, $\tau=0.4$, $B_0=I$, $B_1=0.35I$ and $\mu=\tau^{-1}\mathbf1_{[0,\tau]}\,\dd s$. We use $N$ equal-mass midpoint atoms, $\mathbb P^1$ finite elements, backward Euler and the piecewise affine time reconstruction. Every rectangle is split into two right triangles, with the top-left and bottom-right diagonals reversed to avoid a triangle having three vertices on the boundary (\Cref{fig:corner-mesh}); the mesh of \Cref{prop:eq-obstruction} omits this correction. Here, $d_{\mu,0}=0.65$ and
\begin{align}\label{eq:2d-transport-quantities}
W_1(\mu,\mu_N)=\frac{\tau}{4N},\quad
\|\Delta_{\mu_N,\mu}\|_{L^2(0,\Tend)}=\frac1N\sqrt{\frac{\tau}{12}}.
\end{align}
These values enter the acceptance test \eqref{eq:quantile-acceptance-test}. The manufactured solution and forcing are
\begin{align}\label{eq:2d-manufactured}
 u_\mu(x,y,t)&=e^t\sin(\pi x)\sin(\pi y),\\
 f(x,y,t)&=e^t\sin(\pi x)\sin(\pi y)
 \left[1+2\pi^2\left(1+0.35\frac{1-e^{-\tau}}{\tau}\right)\right].\notag
\end{align}
The history is $\psi(t)=e^t\sin(\pi x)\sin(\pi y)$ for $-\tau<t<0$ and $u_0=\psi(0)$, so that $u_\mu$ is the exact solution; this history produces the factor $(1-e^{-\tau})/\tau$ in $f$. The initial numerical state is the nodal interpolant of $u_0$, so the history junction need not be continuous. The finite element load and the mixed RT reconstruction are assembled with the seven-point degree-five triangle rule $Q_T^{(7)}$. At each time quadrature point, we form
\begin{align*}
 \hat g_h(t)|_T:=\frac{1}{|T|_d}Q_T^{(7)}[g(t)]
\end{align*}
and compute a global $\mathbb{RT}^{0}$ flux by the constrained least-squares system
\begin{align*}
 \min_{\substack{\sigma_h\in\mathbb{RT}^{0}(\mathbb T_h)\\
                  \operatorname{div}\sigma_h=\hat g_h(t)}}
 \|F_N(t)-\sigma_h\|_{L^2(\Omega)^2},
\end{align*}
with its integrals assembled by the same seven-point rule. To evaluate the estimator, we approximate the exact elementwise means $\Pi_h^0g$, the oscillation norm and the flux mismatch by a $14\times14$ Gauss--Legendre tensor rule mapped to each triangle. In time, a six-point Gauss rule is applied separately on intervals cut at time-grid nodes and at the shifted atomic delay nodes. The divergence defect is obtained from the assembled RT divergence matrix after the mixed solve. We use $C_{P,\Omega}=1/(\pi\sqrt2)$, which is valid on the unit square by \Cref{rem:explicit-cp}. The reported $E$ is a quadrature approximation of $\|u_\mu-u_h\|_{L^2(0,\Tend;V)}$; its evaluation uses the seven-point triangle and three-point time rules.

With superscript $(Q)$ denoting numerical evaluation, we report the initial-data term and all five residual contributions in \eqref{eq:computable-majorant}:
\begin{align*}
 B^{(Q)}:=\frac{\eta_{\mathrm{init}}^{(Q)}}{\sqrt{2d_{\mu,0}}}
 +\frac{\eta_{\mathrm{flux}}^{(Q)}+\eta_{\mathrm{osc}}^{(Q)}
 +\eta_{\mathrm{mean}}^{(Q)}+\hat\eta_{\mathrm{div}}^{(Q)}
 +\eta_{\mathrm{ker}}^{\mathrm{cdf},(Q)}}{d_{\mu,0}}.
\end{align*}
In particular, $\eta_{\mathrm{osc}}^{(Q)}$ uses the numerically evaluated $L^2$ projection $\Pi_h^0g$, not $\hat g_h$. Repeating the calculations with a $10\times10$ spatial rule and a four-point time rule agrees in the displayed digits. For the manufactured tests, the evaluated $\eta_{\mathrm{mean}}$ ranges from $1.36\times10^{-9}$ to $8.64\times10^{-8}$, while $\hat\eta_{\mathrm{div}}<1.2\times10^{-14}$. Although every contribution is included, $B^{(Q)}$ is not a verified numerical enclosure of the exact right-hand side.

\begin{table}[htbp]
\centering
\caption{Manufactured solution: approximate error $E$, quadrature evaluation $B^{(Q)}$ of the full majorant (not a validated numerical enclosure), and $\eta_K:=\eta_{\mathrm{ker}}^{\mathrm{cdf},(Q)}$.}
\label{tab:manufactured-errors}
\setlength{\tabcolsep}{4pt}
\begin{tabular}{@{}rrrrrrrr@{}}\hline
$n_x$ & $n_y$ & $L$ & $N$ & $E$ & $B^{(Q)}$ & $B^{(Q)}/E$ & $\eta_K$\\\hline
8&8&20&8&0.60763&1.59175&2.620&0.01740\\
12&12&20&8&0.41065&1.03455&2.519&0.01642\\
16&16&20&8&0.31118&0.77033&2.476&0.01590\\
12&12&40&8&0.40766&1.03147&2.530&0.01659\\
12&12&20&4&0.41062&1.05985&2.581&0.03284\\
12&12&20&16&0.41066&1.02190&2.488&0.00821\\
\hline
\end{tabular}
\end{table}
At $12\times12$, the delay term decreases from $0.03284$ to $0.00821$, in proportion to $N^{-1}$, as $N$ increases from $4$ to $16$, whereas the error stays near $0.411$.  The cancellation-preserving residual $\eta_{\mathrm{ker}}^{\mathrm{rec}}$ is $0.002316$, compared with the cumulative-mass value $0.016418$, for $N=8$ on this mesh.

To evaluate the same complete majorant on a problem without a known solution, we also take $u_0=0$, $\psi=0$ and
\begin{align}\label{eq:2d-local-heat}
f(x,y,t)=(1-e^{-8t})\exp[-80((x-0.30)^2+(y-0.65)^2)].
\end{align}
Here, $J_h=\eta_{\mathrm{init}}=0$. For $N=8$, $L=20$ the quadrature evaluations $B^{(Q)}$ of the complete majorant on $8\times8$, $12\times12$ and $16\times16$ are respectively $0.011162$, $0.007183$ and $0.005298$; the associated $\eta_K$ values are $0.000113$, $0.000116$ and $0.000118$. On the $8\times8$ mesh, $\eta_{\mathrm{ker}}^{\mathrm{rec}}=0.00000483$, compared with $\eta_K=0.00011339$. The reference error is unknown.

\subsection{Grid-aligned RT fluxes}
\label{subsec:aligned-numerics}
The efficiency argument assumes that the delay atoms lie on the time grid. We use the fully coupled refinement $(n_x,n_y,L,N)=(p,p^2,4p,p)$, where $k=0.2/p$ and the midpoint delays satisfy $s_j=(2j-1)k$. Both spatial spacings tend to zero while $h_x/h_y=p\to\infty$; the time step and atomic spacing also tend to zero. The RT divergence data use the time average $f_k$. In the computation, the elementwise projection is approximated by the seven-point triangle rule,
\begin{align}\label{eq:aligned-qhat}
 \hat q_h^n|_T:=\frac{1}{|T|_d}Q_T^{(7)}[g^n],
 \quad q^n:=\Pi_h^0g^n.
\end{align}
The stepwise RT flux $\sigma_h|_{I_n}$ is obtained from the constrained global least-squares solve with divergence data $\hat q_h^n$. The equivalent Marini CR post-processing with the same divergence data gives edge degrees of freedom differing by less than $9\times10^{-9}$ (\Cref{lem:crrt-optimality}); the theoretical conditional-efficiency result of \Cref{thm:conditional-efficiency} concerns the exact data $q^n$.
For each $\beta\ge0$, put
\begin{align}\label{eq:aligned-defect-terms}
 \eta_{\mathrm{mean},k}(\beta)
 &:=C_P\|q^n-\hat q_h^n\|_{L^2_\beta(0,\Tend;L^2(\Omega))},\\
 \hat\eta_{\mathrm{div},k}(\beta)
 &:=C_P\|\hat q_h^n-\operatorname{div}\sigma_h\|_{L^2_\beta(0,\Tend;L^2(\Omega))},
\end{align}
where the quantities inside the norms are constant on each $I_n$ and $C_P=1/\pi$. With $\beta=0$, we define the complete stepwise majorant by
\begin{align}\label{eq:Bstep}
 B_{\rm step}:=\frac{\eta_{\rm init}}{\sqrt{2d_{\mu,0}}}
 +\frac{\eta_{\rm flux}+\eta_{\mathrm{osc},k}+\eta_{f,k}
 +\eta_{\mathrm{mean},k}+\hat\eta_{\mathrm{div},k}
 +\eta_{\rm ker}^{\rm rec}}{d_{\mu,0}}.
\end{align}
Here, $\eta_{\rm flux}=\eta_{\rm flux}(\sigma_h)$ uses the computed stepwise flux, $\eta_{\mathrm{osc},k}$ is formed using the true element mean $q^n$ (approximated by higher-order quadrature for the tables), and $\eta_{f,k}$ is as in \Cref{sec:efficiency}. The term $\eta_{\rm ker}^{\rm rec}$ is the cancellation-preserving residual of \Cref{prop:cancellation-residual}. To justify \eqref{eq:Bstep}, decompose the step residual into the flux mismatch, the elementwise oscillation $g^n-q^n$, the mean discrepancy $q^n-\hat q_h^n$, the divergence defect $\hat q_h^n-\operatorname{div}\sigma_h$, and the temporal forcing discrepancy $f^n-f(t)$. The elementwise and global Poincar\'e inequalities then give
\begin{align*}
 \|R_{\mu_N}\|_{L^2_\beta(0,\Tend;V^*)}
 \le\eta_{\rm flux}+\eta_{\mathrm{osc},k}+\eta_{f,k}
 +\eta_{\mathrm{mean},k}+\hat\eta_{\mathrm{div},k}.
\end{align*}
Combining this estimate with \Cref{thm:fem-certificate,prop:cancellation-residual} proves that $B_{\rm step}$, with exact norms and integrals, bounds $\|u_\mu-u_h\|_{L^2(0,\Tend;V)}$. The table reports its floating-point quadrature evaluation $B_{\rm step}^{(Q)}$, not a validated numerical enclosure. We use the same majorant, comprising the initial contribution and all six residual terms, with $L^2_\beta$ norms and $d_{\mu,\beta}$ for the supercritical tests.
\begin{table}[htbp]
\centering
\caption{Grid-aligned tests ($M=1$) under the joint refinement
$(n_x,n_y,L,N)=(p,p^2,4p,p)$; $B_{\rm step}^{(Q)}$ evaluates all terms in \eqref{eq:Bstep}.}
\label{tab:aligned-subcritical}
\setlength{\tabcolsep}{3pt}
\begin{tabular}{@{}lrrrrrr@{}}\hline
Mesh & $L$ & $N$ & $E$ & $B_{\rm step}^{(Q)}$ & $B_{\rm step}^{(Q)}/E$ & $\eta_{\rm ker}^{\rm rec}$\\\hline
$4\times16$&16&4&0.880099&2.541604&2.89&0.010036\\
$6\times36$&24&6&0.583732&1.588495&2.72&0.004422\\
$8\times64$&32&8&0.436704&1.152934&2.64&0.002478\\
$12\times144$&48&12&0.290580&0.743630&2.56&0.001098\\
$16\times256$&64&16&0.217789&0.548512&2.52&0.000617\\
\hline
\end{tabular}
\end{table}

\subsection{Supercritical mass and the limits of efficiency}
\label{subsec:supercritical-2d}
To test \Cref{cor:weight-optimised-certificate}, we use
\begin{align}\label{eq:supercritical-2d-law}
\dd\mu(s)=\frac{M}{0.2}\mathbf1_{[0.2,0.4]}(s)\,\dd s,
\quad M=4.
\end{align}
With the same history, we choose the forcing so that the exact solution for the original law remains $e^t\sin(\pi x)\sin(\pi y)$; its coefficient is obtained from \eqref{eq:2d-manufactured} by replacing $(1-e^{-\tau})/\tau$ with $\frac M{0.2}\int_{0.2}^{0.4}e^{-s}\,\dd s$. The unweighted margin is $d_{\mu,0}=1-0.35M=-0.4$, yet a positive weighted margin exists. The tests below refine only the mesh, with $N=4$ and $L=32$ fixed; the atoms then lie on the time grid. For each test we select an admissible weight $\beta_{\rm opt}$ numerically and report $B_{\rm step,opt}^{(Q)}$, the quadrature evaluation of $e^{\beta_{\rm opt}\Tend}B_{\rm step}(\beta_{\rm opt})$, where $B_{\rm step}(\beta)$ contains all terms of \eqref{eq:Bstep}, evaluated with $L^2_\beta$ norms and $d_{\mu,\beta}$ in place of $d_{\mu,0}$. The numerical values are computed at the three-decimal weights displayed in the table. The chosen weight need not be a certified global minimiser: by the weighted energy estimate, every admissible weight yields an upper bound for $\|u_\mu-u_h\|_{L^2(0,\Tend;V)}$ when the estimator terms are evaluated rigorously.
\begin{table}[htbp]
\centering
\caption{Grid-aligned supercritical tests ($M=4$) with the
unweighted bound obtained from the weighted majorant $B_{\rm step,opt}^{(Q)}$ at the displayed admissible weights.}
\label{tab:aligned-supercritical}
\setlength{\tabcolsep}{4pt}
\begin{tabular}{@{}lrrrrrr@{}}\hline
Mesh & $L$ & $N$ & $\beta_{\rm opt}$ & $E$ & $B_{\rm step,opt}^{(Q)}$ & $B_{\rm step,opt}^{(Q)}/E$\\\hline
$4\times16$&32&4&2.653&0.878622&18.9495&21.57\\
$8\times64$&32&4&2.671&0.437016&8.7130&19.94\\
$12\times144$&32&4&2.674&0.292276&5.8106&19.88\\
$16\times256$&32&4&2.672&0.220665&4.4571&20.20\\
\hline
\end{tabular}
\end{table}
On the stretched meshes, the reported unweighted majorant remains finite but is conservative, with $B_{\rm step,opt}^{(Q)}/E$ close to $20$. 

Finally, we examine the elliptic (EQ) hypothesis from \Cref{thm:conditional-efficiency}. For $-\Delta z=f$ with the $\mathbb P^1$ Galerkin approximation $z_h$, we write $E_{\rm sp}=\|\nabla(z-z_h)\|_{L^2(\Omega)^2}$, $\eta_{\rm sp}:=\min_{\sigma_h\in\mathcal Q_h(-f)}\|\nabla z_h-\sigma_h\|_{L^2(\Omega)^2}$ and $\mathrm{osc}:=\mathrm{osc}_h(f)$. This is \eqref{eq:elliptic-EQ} with $F_h=\nabla z_h$ and $q=-f$; since $E_{\rm sp}=\|\ell\|_{V^*}$ by \Cref{rem:eq-shapereg}, the last column of \Cref{tab:eq-elliptic-diagnostic} is the ratio in \eqref{eq:elliptic-EQ}. The flux obtained by Marini's post-processing agrees numerically with the constrained global minimum.
\begin{table}[htbp]
\centering
\caption{Elliptic diagnostic: ratios of the minimal RT flux mismatch
to energy error and to error plus data oscillation.}
\label{tab:eq-elliptic-diagnostic}
\setlength{\tabcolsep}{4pt}
\begin{tabular}{@{}llrr@{}}\hline
Solution & Mesh & $\eta_{\rm sp}/E_{\rm sp}$ &
$\eta_{\rm sp}/(E_{\rm sp}+\mathrm{osc})$\\\hline
Smooth&$4\times16$&1.137&0.913\\
Smooth&$8\times64$&1.150&1.028\\
Smooth&$16\times256$&1.154&1.089\\
$y$-layer, $\varepsilon=0.035$&$4\times16$&3.053&1.035\\
$y$-layer, $\varepsilon=0.035$&$8\times64$&3.776&2.099\\
$y$-layer, $\varepsilon=0.035$&$12\times144$&4.164&2.789\\
$y$-layer, $\varepsilon=0.035$&$16\times256$&4.358&3.215\\
$y$-layer, $\varepsilon=0.008$&$8\times64$&5.954&1.151\\
$y$-layer, $\varepsilon=0.008$&$12\times144$&8.035&2.199\\
$y$-layer, $\varepsilon=0.008$&$16\times256$&9.960&3.405\\
\hline
\end{tabular}
\end{table}
The exact solutions are $z=\sin(\pi x)\sin(\pi y)$ and $z=x(1-x)(1-e^{-y/\varepsilon})(1-y)$, respectively. For the layer solutions, the ratios grow under joint refinement and are larger for the thinner layer. These are finite-sample diagnostics; \Cref{prop:eq-obstruction} proves that no aspect-ratio-uniform version of (EQ) holds for the RT flux space.

\FloatBarrier

\section{Conclusion}
\label{sec:conclusion}
We have established computable a posteriori finite element error bounds for diffusion equations with measure-valued delays in the principal part. Although computations use finitely many delay atoms, the bounds concern the solution for the original law. The reliability constants are independent of element aspect ratios and the number of atoms, and require no additional time regularity of the exact solution.

The delay-law approximation is included explicitly in the numerical error budget. Its residual is estimated on the computed time reconstruction, while the cumulative-mass representation retains cancellation and accounts for a jump at the history junction. For quantile atoms, this gives an explicit estimate of the delay-law contribution for a given reconstruction. Exponential weights also yield finite-horizon bounds beyond the mass-subcritical regime, although the constants may deteriorate as delays approach zero.

Under an elliptic equilibration hypothesis, we obtain a two-sided bound for a combined graph error. Our counterexample shows that the RT estimator with prescribed divergence need not be uniformly efficient on anisotropic meshes, even when the flux minimises the least-squares mismatch. Thus, reliability alone gives no uniform efficiency conclusion.

A practical consequence is that the choice of delay atoms can be assessed through its contribution to the error bound, rather than being treated as an unchecked approximation. This provides a basis for developing adaptive methods in which the spatial mesh, time discretisation and atomic approximation are adjusted according to their estimated errors. Such error assessment may also be useful in automated numerical computations. Establishing convergence of these adaptive methods, constructing efficient estimators on general anisotropic meshes, and extending the analysis to partially known delay laws remain subjects for further research.

\section*{Declarations}

\subsection*{Funding}
The author declares that no funds, grants, or other support were received during the preparation of this manuscript.

\subsection*{Data availability statement}
The Python code used to produce the numerical results is available from the author upon reasonable request.

\subsection*{Disclosure statement}
The author declares no competing interests.

\section*{Acknowledgements}
The author used ChatGPT (OpenAI) to assist with language and LaTeX editing, reproducibility code and selected checks of proofs and calculations. The author reviewed the outputs and takes responsibility for the mathematical arguments, numerical results and final manuscript.

\end{document}